\documentclass[11pt]{article}

\usepackage{amsmath, amsthm}
\usepackage{amsmath, amsfonts}
\usepackage{amsmath, amssymb}
\usepackage{amsmath}
\usepackage[all]{xy}

\newtheorem{theorem}{Theorem}[subsection]
\newtheorem{definition}[theorem]{Definition}
\newtheorem{definition-lemma}[theorem]{Definition/Lemma}
\newtheorem{definition-explanation}[theorem]{Definition/Explanation}
\newtheorem{explanation-definition}[theorem]{Explanation/Definition}
\newtheorem{lemma}[theorem]{Lemma}
\newtheorem{lemma-definition}[theorem]{Lemma/Definition}
\newtheorem{proposition}[theorem]{Proposition}
\newtheorem{corollary}[theorem]{Corollary}
\newtheorem{remark}[theorem]{\it Remark}

\newtheorem{example}[theorem]{Example}

\newtheorem{notation}[theorem]{Notation}

\numberwithin{equation}{subsection}

\newtheorem{sdefinition-lemma}[stheorem]{Definition/Lemma}
\newtheorem{sdefinition-explanation}[stheorem]{Definition/Explanation}
\newtheorem{sexplanation-definition}[stheorem]{Explanation/Definition}

\newtheorem{slemma-definition}[stheorem]{Lemma/Definition}

\newcommand{\Aut}{\mbox{\it Aut}\,}

 \newcommand{\Azscriptsize}{\mbox{\scriptsize\it A$\!$z}}

\newcommand{\Br}{\mbox{\it Br}\,}

 \newcommand{\CohCategory}{\mbox{\it ${\cal C}$\!oh}\,}

\newcommand{\Dersheaf}{\mbox{\it ${\cal D}$er}\,}

\newcommand{\End}{\mbox{\it End}\,}
 
\newcommand{\Endsheaf}{\mbox{\it ${\cal E}\!$nd}\,}
 
\newcommand{\Ext}{\mbox{\rm Ext}}
 \newcommand{\footnotesizeExt}{\mbox{\footnotesize\rm Ext}}

\newcommand{\GL}{\mbox{\it GL}}

\newcommand{\Hom}{\mbox{\it Hom}\,}
\newcommand{\Homsheaf}{\mbox{\it ${\cal H}$om}\,}

\newcommand{\Image}{\mbox{\it Im}\,}

\newcommand{\Isom}{\mbox{\it Isom}\,}
\newcommand{\Ker}{\mbox{\it Ker}\,}

\newcommand{\ModCategory}{\mbox{\it ${\cal M}$\!od}\,}
\newcommand{\Mor}{\mbox{\it Mor}\,}

\newcommand{\MorphismCategory}{\mbox{\it ${\cal M}$\!orphism}\,}

\newcommand{\PGL}{\mbox{\it PGL}\,}

\newcommand{\Quot}{\mbox{\it Quot}}

\newcommand{\Scheme}{\mbox{\it ${\cal S}\!$cheme}\,}

\newcommand{\Space}{\mbox{\it Space}\,}

\newcommand{\Spec}{\mbox{\it Spec}\,}
 \newcommand{\boldSpec}{\mbox{\it\bf Spec}\,}

\newcommand{\Supp}{\mbox{\it Supp}\,}
 
\newcommand{\Sym}{\mbox{\it Sym}}

\newcommand{\coker}{\mbox{\it coker}\:}
\newcommand{\degree}{\mbox{\it deg}\,}

\newcommand{\determinant}{\mbox{\it det}\,}
\newcommand{\dimm}{\mbox{\it dim}\,}

\newcommand{\et}{\mbox{\scriptsize\it \'{e}t}\,}

\newcommand{\scriptsizenoncommutative}{\mbox{\scriptsize\rm
                                             noncommutative}}

\newcommand{\pr}{\mbox{\it pr}}
 \newcommand{\prscriptsize}{\mbox{\scriptsize\it pr}}
\newcommand{\productscriptsize}{\mbox{\scriptsize\it product}\,}

\newcommand{\rank}{\mbox{\it rank}\,}

\newcommand{\trace}{\mbox{\it tr}\,}

\begin{document}

\enlargethispage{23cm}

\begin{titlepage}

$ $

\vspace{-1cm}

\noindent\hspace{-1cm}
\parbox{6cm}{\small July 2009}\
   \hspace{8cm}\
   \parbox[t]{5cm}{yymm.nnnn [math.AG]\\ D(5): $B$-field effect}

\vspace{2cm}

\centerline{\large\bf
 Nontrivial Azumaya noncommutative schemes, morphisms therefrom,}
\vspace{1ex}
\centerline{\large\bf
 and their extension by the sheaf of algebras
 of differential operators$\,$:}
\vspace{1ex}
\centerline{\large\bf
 D-branes in a {\boldmath $B$}-field background
 \`{a} la Polchinski-Grothendieck Ansatz}

\bigskip

\vspace{3em}

\centerline{\large
  Chien-Hao Liu
  \hspace{1ex} and \hspace{1ex}
  Shing-Tung Yau
}

\vspace{4em}

\begin{quotation}
\centerline{\bf Abstract}

\vspace{0.3cm}

\baselineskip 12pt  
{\small
 In this continuation of [L-Y1], [L-L-S-Y], [L-Y2], and [L-Y3]
 (arXiv:0709.1515 [math.AG], arXiv:0809.2121 [math.AG],
  arXiv:0901.0342 [math.AG], arXiv:0907.0268 [math.AG]),
 we study D-branes in a target-space(-time)
  with a fixed $B$-field background $(Y,\alpha_B)$
  along the line of the Polchinski-Grothendieck Ansatz,
  explained in [L-Y1] and further extended in the current work.
 We focus first on the gauge-field-twist effect of $B$-field
  to the Chan-Paton module on D-branes.
 Basic properties of the moduli space of D-branes,
  as morphisms from Azumaya schemes with a twisted fundamental module
  to $(Y,\alpha_B)$, are given.
 For holomorphic D-strings, we prove a valuation-criterion property
  of this moduli space.
 The setting is then extended to take into account also
  the deformation-quantization-type noncommutative geometry effect
  of $B$-field to both the D-brane world-volume and
  the superstring target-space(-time) $Y$.
 This brings the notion of twisted ${\cal D}$-modules
   that are realizable as twisted locally-free coherent modules
          with a flat connection
  into the study.
 We use this to realize the notion of
  both the classical and the quantum spectral covers
  as morphisms from Azumaya schemes with
  a fundamental module (with a flat connection in the latter case)
  in a very special situation.
 The 3rd theme
  (subtitled ``Sharp vs.\ Polchinski-Grothendieck") of Sec.~2.2
  is to be read with the work [Sh3] (arXiv:hep-th/0102197) of Sharp
  while
 Sec.~5.2 (subtitled less appropriately
  ``Dijkgraaf-Holland-Su{\l}kowski-Vafa vs.\ Polchinski-Grothendieck")
  is to be read with the related sections in
   [D-H-S-V] (arXiv:0709.4446 [hep-th])  and
   [D-H-S]   (arXiv:0810.4157 [hep-th])
  of Dijkgraaf, Hollands, Su{\l}kowski, and Vafa.
} 
\end{quotation}

\vspace{2em}

\baselineskip 12pt
{\footnotesize
} 

\bigskip

\noindent {\small MSC number 2010:
 14A22, 81T30; 14F10, 14D23, 81T75.
} 

\bigskip

\baselineskip 10pt
{\scriptsize
\noindent{\bf Acknowledgements.}
 We thank
  Andrew Strominger and Cumrun Vafa
   for influencing our understanding of strings and branes over the years.
 For the current D(5), C.-H.L.\ thanks in addition
  Shiraz Minwalla for the numerous stringy inspirations
   during the brewing years;
  Mihnea Popa
   for discussions on D-branes and spectral covers and
  William Oxbury
   for preprint/communication, spring 2002;
  Andrei C\u{a}ld\u{a}raru
   for communicating his thesis and
  Liang Kong and Eric Sharpe
   for sharing their insights on D-branes,
   fall 2007;
  Cumrun Vafa for
   the lecture and answers to questions on D/I-branes and
   ${\cal D}$-modules, December 2007, and the illuminations
   in the topic course in string theory, spring 2009;
  Si Li and Ruifang Song for the participation of D(2), spring 2008;
  Alessandro Tomasiello
   for discussions on $B$-field and $(p,q)$-strings,
  Clark Barwick for the topic course toward $K$-theory and $G$-theory
   and a discussion on twisted $K$-theory, spring 2009;
 and Ling-Miao Chou for the long-standing moral support.
 The project is supported by NSF grants DMS-9803347 and DMS-0074329.
} 

\end{titlepage}

\newpage
\begin{titlepage}
\baselineskip 14pt

$ $

\vspace{12em} 

\centerline{\small\it
 Chien-Hao Liu dedicates this work to his mother,
  Mrs.~Pin-Fen Lin Liu,}
\centerline{\small\it
 and siblings (resp.\ siblings-in-law)}
\centerline{\small\it
 Hsiu-Chuan Liu (Yi-Shou Chao), Chien-Ying Liu (Tracy Tse Chui),}
\centerline{\small\it
 Alice Hsiu-Hsiang Liu (Tom Cheng-Wei Chang), Hsiu-Luan Liu,}
\centerline{\small\it
 and to the memory of his father, Mr.~Han-Tu Liu$^{\dag}$
                                  [$\;^{\dag}\!$deceased 1986],}
\centerline{\small\it
 at a special family moment.}

%

\end{titlepage}

\newpage
$ $

\vspace{-4em}  

\centerline{\sc
 D-Branes in a $B$-Field Background
 \`{a} la Polchinski-Grothendieck Ansatz}

\vspace{2em}

\baselineskip 14pt  

\begin{flushleft}
{\Large\bf 0. Introduction and outline.}
\end{flushleft}
Since the work [Pol1] of Polchinski, D-branes have become
 a central object of study in superstring theory.
It has also motivated numerous related works on the mathematical side.
(Cf. References of [L-Y1], [L-L-S-Y], [L-Y2], [L-Y3]
     for a brief list relevant to the project.)

\bigskip

\begin{flushleft}
{\bf Azumaya structure on D-branes as an origin of D-brany phenomena.}
\end{flushleft}
The emergence of an Azumaya structure on a D-brane world-volume
 follows directly from a comparison of
  (1) {\it the behavior of the open-string-induced field on the D-brane
       that governs its deformation} and
  (2) {\it Grothendieck's contravariant equivalence of
       algebras and local geometries}
 ([Pol2: vol.~I, Sec.~8.7], [L-Y1: Sec.~2],
    and Sec.~2.1).\footnote{Azumaya
                          structure on D-brane world-volume
                            has already been brought to string theorists'
                            attention in late 1990s, for example,
                           from the effect of the background $B$-field
                            on the open-string target-space-time,
                            as explained in [Ka: Sec.~1.2] of Kapustin.
                          However, it should be noted that
                           {\it the emergence of Azumaya structures
                           on a D-brane world-volume
                           is at a more fundamental level than this.
                          It is purely an open-string-induced effect
                           that is enforced on
                           a coincident D-brane world-volume
                           whether or not there is
                            a supersymmetry in the field theory
                             on the D-brane world-volume
                             (or on the open-string world-sheet
                                  with boundary on this brane-world-volume,
                              or on the space-time)
                            or a $B$-field background on the space-time.}
                          Rather, the latter extra SUSY requirement
                           or $B$-field data comes to constraint
                           the class of Azumaya
                           structures that can occur on the D-brane
                           world-volume.
                          For example,
                           in complex geometry language
                            we are discussing
                             holomorphic Azumaya algebra
                              over a complex manifold.
                          This comes from a supersymmetry constraint.
                          $B$-field will then select further
                           a class of such Azumaya algebras,
                           cf.\ [Ka] and Sec.~2.}
 If we turn the history around to take such a structure
  as a fundamental definition of D-branes,  and
  consider morphisms from such objects to a string-target-space(-time),
 then we see that basic D-brane phenomena, e.g.,
   Higgsing/un-Higgsing of gauge field theory on D-brane world-volume,
   deformation and resolution of a singular Calabi-Yau space
   via a D-brane probe,
  can be reproduced;
 [L-Y1], [L-L-S-Y], [L-Y2], and [L-Y3].
 This is an indication that
  Azumaya structure is fundamentally/solidly carved
       into a D-brane world-volume
       as part of its substantial building structures.
 This gives an Azumaya origin of many D-brany phenomena.

\bigskip

\begin{flushleft}
{\bf {\boldmath $B$}-field and
     its effect on fields and geometry in string theory.}
\end{flushleft}
A $B$-field on a space-time $Y$ is a connection on
 a gerbe ${\cal Y}$ over $Y$.
It can be presented
 as a \v{C}ech $0$-cochain $(B_i)_i$ of local $2$-forms $B_i$
 with respect to a cover ${\cal U}=\{U_i\}_i$ on $Y$
 such that on $U_i\cap U_j$,
$B_i-B_j=d\Lambda_{ij}$ for some real $1$-forms $\Lambda_{ij}$
 that satisfies
 $\Lambda_{ij}+\Lambda_{jk}+\Lambda_{ki}=-\sqrt{-1}d\log\alpha_{ijk}$
  on $U_i\cap U_j \cap U_k$,
 where $(\alpha_{ijk})_{ijk}$ is a \v{C}eck $2$-cocycle of
  $U(1)$-valued functions on $Y$
In the algebro-geometric language,
 $(\alpha_{ijk})_{ijk}$ is given by a presentation of
 an equivalence class
 $\alpha_B\in \check{C}^2_{\et}(Y,{\cal O}_Y^{\ast})$
 of \'{e}tale \v{C}ech $2$-cocycles with values in ${\cal O}_Y^{\ast}$.
Through its coupling to the open-string current
 on an open-string world-sheet
 with boundary on a D-brane world-volume $X\subset Y$,
a background $B$-field on $Y$ induces\footnote{Unfamiliar
                          mathematicians are highly recommended
                          to read [Zw] of Zwiebach
                          for a very down-to-earth explanation of this.}
 a twist to the gauge field $A$
 on the Chan-Paton vector bundle $E$ on $X$
 that renders $E$ itself a twisted vector bundle
  with the twist specified by
  $\alpha_B|_X\in \check{C}^2_{\et}(X,{\cal O}_X^{\ast})$.
(Cf.~[Al], [Br], [Ch], [F-W], [Hi2], [Ka], and [Wi1].)
Furthermore,
 the $2$-point functions on the open-string world-sheet
  with boundary on $X$ indicate that
 the D-brane world-volume is deformed
  to a deformation-quantization type noncommutative geometry
  in a way that is governed by the $B$-field
  (and the space-time metric).
(Cf.\ [C-H1, C-H2], [C-K], [Schmo], and [S-W].)

\bigskip

\noindent
{\it Remark 0.1.\ $[$open-string world-sheet anomaly$]$.}
 There is a further effect on D-branes that arises from
  the global world-sheet anomaly on the open-string world-sheet
  with the boundary on the D-brane world-volume ([F-W]).
 This anomaly effect is ignored in the current work.
 See ibidem, [C-K-S], [Ka], and [K-S] for more discussions.

\bigskip

\begin{flushleft}
{\bf D-brane as a master object in superstring theory
     vs.\ morphism from Azumaya schemes with a fundamental module
          as a master object in geometry.}
\end{flushleft}
In this continuation of [L-Y1], [L-L-S-Y], [L-Y2], and [L-Y3],
we study D-branes in a fixed $B$-field background $(Y,\alpha_B)$
 along the line of the Polchinski-Grothendieck Ansatz,
 explained in [L-Y1] and further extended in the current work.
We focus first on the twist effect of $B$-field to
 the Chan-Paton module on D-branes.
Basic properties of the moduli space of D-branes,
 as morphisms from Azumaya schemes with a twisted fundamental module
 to $(Y,\alpha_B)$, are given.
For holomorphic D-strings, we prove a valuation-criterion property
 of this moduli space.
The setting is then extended to take into account also
 the deformation-quantization effect of $B$-field to
 both the D-brane world-volume and the target-space $Y$.
This brings the notion of twisted ${\cal D}$-modules
  that are realizable as twisted locally-free coherent modules
   with a flat connection
 into the study.
We use this to realize the notion of
 both the classical and the quantum spectral covers
 as morphisms from Azumaya schemes with a fundamental module
  (with a flat connection in the latter case)
 in a very special situation.
The 3rd theme
 (subtitled ``{\it Sharp vs.\ Polchinski-Grothendieck}") of Sec.~2.2
 is to be read with the work [Sh3] of Sharp
 while
Sec.~5.2 (subtitled less appropriately
 ``{\it Dijkgraaf-Holland-Su{\l}kowski-Vafa
        vs.\ Polchinski-Grothendieck}")
 is to be read with the related sections in
  [D-H-S-V]  and
  [D-H-S]
 of Dijkgraaf, Hollands, Su{\l}kowski, and Vafa.
{From} this, we see once again:
 \begin{itemize}
  \item[$\cdot$] {\it
   the master nature of
     morphisms from Azumaya schemes with a fundamental module\\
    in geometry
   in parallel to the master nature of D-branes in superstring theory.}
 \end{itemize}
This would be highly surprising/un-anticipated on the mathematics side
 if not because of the Polchinski-Grothendieck Ansatz,
 which realizes morphisms from Azumaya manifolds/schemes/ stacks
  with a fundamental module
  as the {\it lowest level} presentation of D-branes,
 and superstring theory dictates the master nature of such an object.
Together with [L-L-S-Y] (D(2)), [L-Y2] (D(3)), and [L-Y3] (D(4)),
 the following diagram of unity emerges:

 \medskip

 $$
  
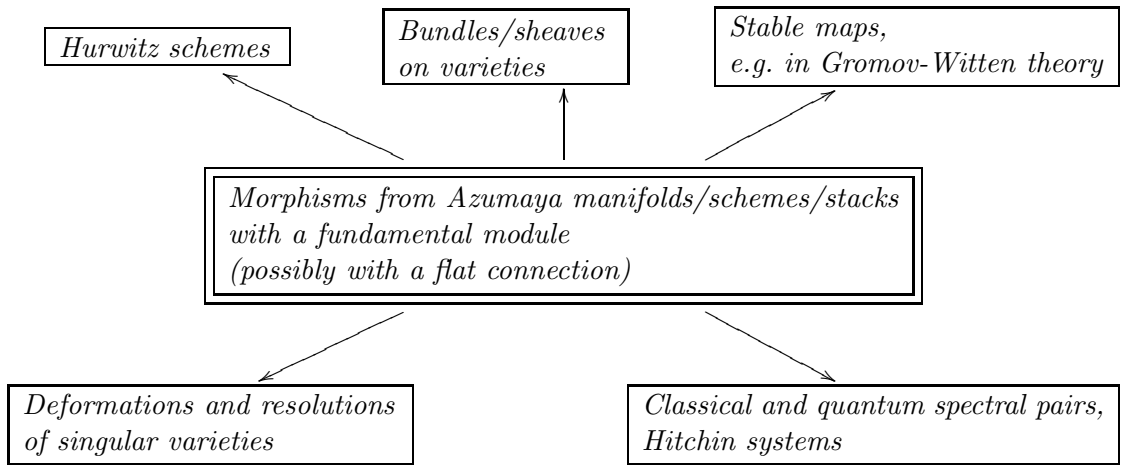
\hspace{-3em}
  \xymatrix @R=6ex @C=-18ex {
   \framebox[8.5em][c]{\parbox{7.5em}{\it Hurwitz schemes}}
    & \hspace{-4em}
      \framebox[8.5em][c]{\parbox{7.5em}{\it
       Bundles/sheaves\\ on varieties}}
    & \framebox[14em][c]{\parbox{13em}{\it
       Stable maps,\\ e.g.\ in Gromov-Witten theory}}        \\
   & \framebox[24.8em]{\framebox[24.2em][c]{\parbox{23.2em}{\it
         Morphisms
         from Azumaya manifolds/schemes/stacks\\
         with a fundamental module\\
         (possibly with a flat connection)}}}
     \ar[lu] \ar[u] \ar[ru] \ar[ld] \ar[rd]                  \\
   \hspace{3em}\framebox[14em][c]{\parbox{13em}{\it
    Deformations and resolutions\\ of singular varieties}}
    && \hspace{-3em}
       \framebox[17em][c]{\parbox{16em}{\it
        Classical and quantum spectral pairs,\\ Hitchin systems}}
  }
 $$

 \bigskip

\noindent It is anticipated that this is only a part of a
to-be-understood
  complete diagram of unity in geometry
 in view of the ubiquity of D-branes in superstring theory.

\bigskip
\bigskip

\noindent
{\bf Convention.}
 Standard notations, terminology, operations, facts in
  (1) physics aspects of D-branes ;
  (2) algebraic geometry and stacks
  can be found respectively in
  (1) [Pol2], [Jo], and [Zw];
  (2) [Ha] and [L-MB].
 \begin{itemize}
  \item[$\cdot$]
   All schemes are Noetherian over ${\Bbb C}$
   unless otherwise noted..

  \item[$\cdot$]
   $B$-{\it field} (in the sense of quantum field theory)
   vs.\ {\it base} scheme $B$ vs.\ D-branes of {\it type} B.

  \item[$\cdot$]
   D-branes for {\it Dirichlet}-branes
    vs.\ ${\cal D}$-modules for modules
         of the sheaf ${\cal D}$
         of algebras of {\it differential} operators.

  \item[$\cdot$]
   The word ``{\it twist/twisting}" has two different meanings:
    (1) in the sense of twisted sheaves as a presentation of
       sheaves on gerbes and
    (2) the operation of tensoring by (usually)
        a (twisted or ordinary in the sense of (1)) line bundle.
 \end{itemize}

\bigskip

\begin{flushleft}
{\bf Outline.}
\end{flushleft}
{\small
\baselineskip 12pt  
\begin{itemize}
 \item[0.]
  Introduction.
  \vspace{-.6ex}
  \begin{itemize}
   \item[$\cdot$]
    Azumaya structure on D-branes as an origin of D-brany phenomena.

   \item[$\cdot$]
    $B$-field and
    its effect on fields and geometry in string theory

   \item[$\cdot$]
    D-brane as a master object in superstring theory
     vs.\ morphism from Azumaya schemes with a fundamental module
          as a master object in geometry.
  \end{itemize}

 \item[1.]
  Gerbes, twisted sheaves, and Azumaya algebras over a scheme.
  \vspace{-.6ex}
  \begin{itemize}
   \item[1.1]
    Gerbes and twisted sheaves over a scheme.

   \item[1.2]
    (General) Azumaya algebras over a scheme.
  \end{itemize}

 \item[2.]
  Azumaya geometry and D-branes \`{a} la Polchinski-Grothendieck Ansatz
   revisited:\\ the twist from a $B$-field background.
  \vspace{-.6ex}
  \begin{itemize}
   \item[2.1]
    Polchinski-Grothendieck Ansatz revisited with the \'{e}tale topology.

   \item[2.2]
    D-branes in a $B$-field background as morphisms
    from Azumaya schemes with a twisted\\ fundamental module.
  \end{itemize}

 \item[3.]
 The moduli stack of morphisms.
  \vspace{-.6ex}
  \begin{itemize}
   \item[3.1]
    Family of D-branes in a $B$-field background,
    twisted Hilbert polynomials, and boundedness.

   \item[3.2]
    ${\frak M}_{A\!z(X_S/S,\alpha_S)^f}(Y,\alpha_B)$ is algebraic.
  \end{itemize}

 \item[4.]
 The case of holomorphic D-strings.
  \vspace{-.6ex}
  \begin{itemize}
   \item[4.1]
    The moduli stack ${\frak M}_{A\!z(g,r,\chi)^f}(Y,\alpha_B;\beta)$
    of morphisms from Azumaya prestable curves\\ to $(Y,\alpha_B)$.

   \item[4.2]
    Fillability/valuation-criterion property of
    ${\frak M}_{A\!z(g,r,\chi)^f}(Y,\alpha_B;\beta)$.
  \end{itemize}

 \item[5.]
 The extension by the sheaf ${\cal D}$ of differential operators.
  \vspace{-.6ex}
  \begin{itemize}
   \item[5.1]
    Azumaya schemes with a fundamental module
    with a flat connection.

   \item[5.2]
    Deformation quantizations of spectral covers
    in a cotangent bundle.
  \end{itemize}
\end{itemize}
} 

\newpage

\section{Gerbes, twisted sheaves, and Azumaya algebras over a scheme.}
To fix terminology and notations,
essential definitions of gerbes and twisted sheaves
 are given in this section.
Readers are referred to
 [Br: Chap.~5], [C\u{a}: Chap.~1],
  [Lie1: Chap.~2], [Mi: Chap.~IV]
 and also [Ch], [Gir], [Hi2] for further details  and
 to, e.g., [Sh2] and [C-K-S] to get a glimpse
  of gerbes and twisted sheaves in string theory.
We will assume that the schemes $X$ and $Y$ in the following discussions
 are quasi-projective\footnote{This
                               technical assumption is imposed to render
                               the \'{e}tale \v{C}ech cohomology
                               $\check{H}^{\ast}_{et}(X,{\cal F})$
                               and the \'{e}tale cohomology
                               $H_{et}^{\ast}(X,{\cal F})$ identical
                               for ${\cal F}$ a sheaf on $X_{\et}$.}
 over ${\Bbb C}$.

\bigskip

\subsection{Gerbes and twisted sheaves over a scheme.}

\begin{flushleft}
{\bf Gerbes over a scheme and coherent sheaves thereupon.}
\end{flushleft}
Let $X$ be a (Noetherian) scheme (over ${\Bbb C}$).
Given a stack ${\cal S}$
  over the category $\Scheme/X$ of schemes over $X$,
 we will denote the groupoid ${\cal S}(U)$ assigned by ${\cal S}$ to
 a $(U\rightarrow X)\in \Scheme/X$ also by ${\cal S}_U$.
An element $s\in {\cal S}_U$ will be called a {\it section} of
 ${\cal S}$ over $U$.
$s$ defines a morphism $s:U\rightarrow {\cal S}$, and conversely.
Thus, we will denote $s\in {\cal S}_U$ and $s:U\rightarrow {\cal S}$
 interchangeably.
We will equip $\Scheme/X$ with the fppf topology
 unless otherwise noted.
This induces a topology on a stack over $\Scheme/X$.

\smallskip

\begin{definition}
{\bf [gerbe over $X$].} {\rm
 A {\it gerbe over $X$} is a stack ${\cal X}$ over $\Scheme/X$
  that has the following two properties:
  \begin{itemize}
   \item[(1)]
   {\it \'{e}tale local existence of a section}$\,$:
    For any $U\rightarrow X$,
     there exists an \'{e}tale cover $U^{\prime}\rightarrow U$ of $U$
     such that ${\cal X}_{U^{\prime}}$ is nonempty.

   \item[(2)]
   {\it sections \'etale locally isomorphic}$\,$:
    For any $U\rightarrow X$ and $s_1$, $s_2\in {\cal X}_U$,
     there exists an \'{e}tale cover $p: U^{\prime}\rightarrow U$ of $U$
     such that $p^{\ast}s_1\simeq p^{\ast}s_2$ in ${\cal X}_{U^{\prime}}$.
  \end{itemize}
 We will denote a gerbe ${\cal X}$ over $X$ also by ${\cal X}/X$
   to manifest the {\it underlying scheme}\footnote{The
                                reason we call $X$
                                 the {\it underlying scheme} of
                                 the gerbe ${\cal X}$ is that
                                when ${\cal X}$ arises
                                  as the moduli stack
                                   of a moduli problem of
                                   a class of objects,
                                 the scheme $X$ becomes
                                  the coarse moduli space of
                                  the moduli problem.
                                $X$ parameterizes all the ${\Bbb C}$-points
                                 of ${\cal X}$
                                while ${\cal X}$ encodes in addition
                                 the data of automorphisms of the objects
                                 these ${\Bbb C}$-points represent.
                                See Definition/Lemma~1.2.3 for an example.}
                                %
                               $X$,
  particularly when there are different underlying schemes
   involved in the discussion.
}\end{definition}

\smallskip

\begin{lemma-definition}
{\bf [sheaf of automorphism groups and its right action].}
{\rm
 Given a gerbe ${\cal X}/X$,
  the assignment
   $s\in {\cal X}_U \mapsto \Aut(s) =\Mor(s,s)\subset \Mor({\cal X}_U)$
  is a sheaf ${\cal A}({\cal X})$ on the stack ${\cal X}$.
 We will call ${\cal A}({\cal X})$ the {\it sheaf of automorphism groups}
  on ${\cal X}$.
 Let ${\cal F}$ be a sheaf on ${\cal X}$.
 Then the operation of pulling-back by automorphisms
  defines a natural right group action\footnote{Explicitly,
                             let $s\in {\cal X}_U$, $f\in {\cal F}(s)$,
                              and $h\in\Aut(s)$.
                             Then, $\mu(f,h)= h^{\ast}f\in {\cal F}(s)$.}
  $\mu:{\cal F}\times {\cal A}({\cal X})\rightarrow {\cal F}$
  of ${\cal A}({\cal X})$ on ${\cal F}$.
}\end{lemma-definition}

\smallskip

Let
 ${\cal X}$ be a gerbe over $X$ and
 ${\cal O}_X^{\ast}$ be the sheaf of invertible elements
  of ${\cal O}_X$.
Denote the pull-back of ${\cal O}_X^{\ast}$ to ${\cal X}$
 via the structure morphism ${\cal X}\rightarrow X$
 also by ${\cal O}_X^{\ast}$.
This is the sheaf on ${\cal X}$
 that assigns to each $s\in {\cal X}_U$
 the (multiplicative) abelian group ${\cal O}_U^{\ast}(U)$.

\smallskip

\begin{definition}
{\bf [gerbe with band ${\cal O}_X^{\ast}$].} {\rm
 A {\it gerbe} over $X$ {\it with band} ${\cal O}_X^{\ast}$
  is a gerbe ${\cal X}/X$ with an isomorphism
  ${\cal O}_X^{\ast}\stackrel{\sim}{\rightarrow} {\cal A}({\cal X})$.
}\end{definition}

\smallskip

\begin{lemma}
{\bf [gerbe as algebraic stack].}
{\rm ([Lie1: Lemma~2.2.1.1].)}
 Let ${\cal X}$ be a gerbe over $X$ with band ${\cal O}_X^{\ast}$.
 Then ${\cal X}$ is an (Noetherian) algebraic stack\footnote{I.e.\
                                                   Artin stack.}
      over $\Scheme/X$.
\end{lemma}

\smallskip

\noindent
An atlas for such an ${\cal X}/X$ is given by an \'{e}tale cover
  $U\rightarrow X$
 with the property that
  for any $x\in X$,
   there is a connected component $U_x$ of $U$ that
    gives an \'{e}tale neighborhood of $x\in X$
    such that ${\cal X}_{U_x}$ is nonempty.

The notion of a (Cartesian) coherent sheaf ${\cal F}$
 on a gerbe ${\cal X}/X$ with band ${\cal O}_X^{\ast}$
 is then defined as that for algebraic stacks, given in [L-MB].
Its {\it (stacky) support} $\Supp{\cal F}$,
 defined by the annihilator ideal sheaf
  $\Ker( {\cal O}_{\cal X} \rightarrow
              \Endsheaf_{{\cal O}_{\cal X}}({\cal F}))$
 of ${\cal F}$, is a closed substack of ${\cal X}$.

\bigskip

\begin{flushleft}
{\bf Twisted sheaves \`{a} la C\u{a}ld\u{a}raru.}
\end{flushleft}
Given an \'{e}tale cover $p:U^{(0)}:=\amalg_{i\in I}U_i\rightarrow X$
  of $X$,
 we will adopt the following notations:
 \begin{itemize}
  \item[$\cdot$]
   $U_{ij}\, :=\, U_i\times_X U_j\, =:\, U_i\cap U_j\,$,
   $\; U_{ijk}\, :=\, U_i\times_X U_j\times_X U_k\,
                  =:\, U_i\cap U_j\cap U_k\,$;

  \item[$\cdot$]
   $\xymatrix{
     \cdots\;\ar@<1.2ex>[r] \ar@<.4ex>[r]
             \ar@<-.4ex>[r] \ar@<-1.2ex>[r]
      & U^{(2)}
        := U\times_X U\times_X U
    }$\\
   $\xymatrix{
      \hspace{10em}
      \ar@<.8ex>[rrr]^-{p_{12},\, p_{13},\, p_{23}}
           \ar[rrr] \ar@<-.8ex>[rrr]
      &&& U^{(1)}
          := U\times_X U
          \ar@<.4ex>[rr]^-{p_1,\, p_2} \ar@<-.4ex>[rr]
      && U^{(0)} \ar[r]^-p & X
     }$\\ \\
    are the projection maps from fibered products as indicated;
   the restriction of these projections maps to respectively
    $U_{ijk}$ and $U_{ij}$ will be denoted the same;

  \item[$\cdot$]
   the pull-back of an ${\cal O}_{U_i}$-module ${\cal F}_i$ on $U_i$
    to $U_{ij}$, $U_{ji}$, $U_{ijk}, \,\cdots\,$
    via compositions of these projection maps will be denoted by
   ${\cal F}_i|_{U_{ij}}$, ${\cal F}_i|_{U_{ji}}$,
   ${\cal F}_i|_{U_{ijk}}$, $\cdots\,$ respectively.
 \end{itemize}

\smallskip

\begin{definition}
{\bf [$\alpha$-twisted ${\cal O}_X$-module on an \'{e}tale cover of $X$].}
{\rm ([C\u{a}: Definition~1.2.1].)} {\rm
 Let $\alpha\in \check{C}^2_{\et}(X,{\cal O}_X^{\ast})$
  be a \v{C}ech $2$-cocycle in the \'{e}tale topology of $X$.
 An {\it $\alpha$-twisted ${\cal O}_X$-module on an \'{e}tale cover of}
  $X$ is a triple
  $$
   {\cal F}\; =\; ( \{U_i\}_{i\in I},\,
                    \{{\cal F}_i\}_{i\in I},\,
                    \{\phi_{ij}\}_{i,j\in I} )
  $$
  that consists of the following data
  \begin{itemize}
   \item[$\cdot$]
    an \'{e}tale cover
     $p:U:=\amalg_{i\in I}U_i\rightarrow X$ of $X$
     on which $\alpha$ can be represented as a $2$-cocycle:\\
     $$
      \alpha\; =\;
       \{\, \alpha_{ijk}\,:\,
        \alpha_{ijk} \in \Gamma( U_{ijk}, {\cal O}_X^{\ast} )\;
       \mbox{with
        $\alpha_{jkl}\alpha_{ikl}^{-1}\alpha_{ijl}\alpha_{ijk}^{-1}=1$
        on $U_{ijkl}$ for all $i,j,k,l\in I$} \,\}\,,
     $$
    such a cover will be called
     an {\it $\alpha$-admissible \'{e}tale cover} of $X$;

   \item[$\cdot$]
    ${\cal F}_i$ is a sheaf of ${\cal O}_{U_i}$-modules on $U_i$;

   \item[$\cdot$] ({\it gluing data})$\;$
    $\phi_{ij}: {\cal F}_i|_{U_{ij}} \rightarrow {\cal F}_j|_{U_{ij}}$
     is an ${\cal O}_{U_{ij}}$-module isomorphism
    that satisfies
    \begin{itemize}
     \item[(1)]
      $\phi_{ii}$ is the identity map for all $i\in I$;

     \item[(2)]
      $\phi_{ij}=\phi_{ji}^{-1}$ for all $i,j\in I$;

     \item[(3)] ({\it twisted cocycle condition})\;
      $\phi_{ki}\circ \phi_{jk}\circ \phi_{ij}$
       is the multiplication by $\alpha_{ijk}$ on ${\cal F}_i|_{U_{ijk}}$.
    \end{itemize}
  \end{itemize}
 ${\cal F}$ is said to be {\it coherent}
   (resp.\ {\it quasi-coherent}, {\it locally free})
  if ${\cal F}_i$ is a coherent
   (resp.\ quasi-coherent, locally free)
   ${\cal O}_{U_i}$-module for all $i\in I$.
 A {\it homomorphism}
   $$
    h\;:\; {\cal F}\;=\; ( \{U_i\}_{i\in I},\,
                     \{{\cal F}_i\}_{i\in I},\,
                     \{\phi_{ij}\}_{i,j\in I} )\;
      \longrightarrow\;
     {\cal F}^{\prime}\; =\; ( \{U_i\}_{i\in I},\,
                     \{{\cal F}^{\prime}_i\}_{i\in I},\,
                     \{\phi^{\prime}_{ij}\}_{i,j\in I} )
    $$
   between $\alpha$-twisted ${\cal O}_X$-modules
   on the \'{e}tale cover $p$ of $X$
  is a collection
   $\{h_i:{\cal F}_i\rightarrow {\cal F}^{\prime}_i\}_{i\in I}$,
    where $h_i$ is an ${\cal O}_{U_i}$-module homomorphism,
   such that $\phi^{\prime}_{ij}\circ h_i= h_j\circ \phi_{ij}$
    for all $i,j\in I$.
 In particular, $h$ is an {\it isomorphism} if all $h_i$ are isomorphisms.
 Denote by $\ModCategory(X,\alpha, p)$
  the category of $\alpha$-twisted ${\cal O}_X$-modules
  on the \'{e}tale cover $p:U^{(0)}\rightarrow X$ of $X$.
}\end{definition}

\smallskip

Given an $\alpha$-twisted sheaf ${\cal F}$ on the \'{e}tale cover
  $p:U\rightarrow X$ of $X$.
let $p^{\prime}:U^{\prime}\rightarrow X$
  be an \'{e}tale refinement of $p:U\rightarrow X$.
Then
 $\alpha$ can be represented also on $p^{\prime}:U^{\prime}\rightarrow X$
  and
 ${\cal F}$ on $p$ defines an $\alpha$-twisted ${\cal O}_X$-module
  ${\cal F}^{\prime}$ on $p^{\prime}$ via the pull-back
 under the built-in \'{e}tale cover $U^{\prime}\rightarrow U$ of $U$.
This defines an equivalence of categories:
 $$
  \ModCategory(X,\alpha,p)\;
   \longrightarrow\; \ModCategory(X,\alpha,p^{\prime})\,.
 $$
([C\u{a}: Lemma~1.2.3, Lemma~1.2.4, Remark~1.2.5].)

\smallskip

\begin{definition}
{\bf [$\alpha$-twisted ${\cal O}_X$-module on $X$].} {\rm
 An {\it $\alpha$-twisted ${\cal O}_X$-module on $X$}
  is an equivalence class $[{\cal F}]$
  of $\alpha$-twisted ${\cal O}_X$-modules ${\cal F}$
  on \'{e}tale covers of $X$,
  where the equivalence relation is generated by \'{e}tale refinements
   and descents by
   \"{e}tale covers of $X$ on which $\alpha$ can be represented.
 An ${\cal F}^{\prime}\in [{\cal F}]$ is called a {\it representative}
  of the $\alpha$-twisted ${\cal O}_X$-module $[{\cal F}]$.
 For simplicity of terminology,
  we will also call ${\cal F}^{\prime}$ directly
  an $\alpha$-twisted ${\cal O}_X$-module on $X$.
}\end{definition}

\smallskip

\noindent
Cf.\ [C\u{a}: Corollary~1.2.6 and Remark~1.2.7].

Standard notions of ${\cal O}_X$-modules,
 in particular
  \begin{itemize}
   \item[$\cdot$]
    the {\it scheme-theoretic support} $\Supp{\cal E}$,

   \item[$\cdot$]
    the {\it dimension} $\dimm{\cal E}$, and

   \item[$\cdot$]
    {\it flatness over a base $S$}
  \end{itemize}
 of an $\alpha$-twisted sheaf ${\cal E}$ on $X$ (or on $X/S$)
 are defined via a(ny) presentation of ${\cal E}$
 on an $\alpha$-admissible \'{e}tale cover $U\rightarrow X$.

Standard operations on ${\cal O}_{\bullet}$-modules apply to
 twisted ${\cal O}_{\bullet}$-modules
 on appropriate admissible \'{e}tale covers
 by applying the operations component by component over the cover.
These operations apply then to twisted ${\cal O}_X$-modules as well:
They are defined on representatives of twisted sheaves
 in such a way that they pass to each other
  by pull-back and descent under \'{e}tale refinements
  of admissible \'{e}tale covers.
In particular:

\smallskip

\begin{proposition}
{\bf [basic operations on twisted sheaves].}
{\rm ([C\u{a}: Proposition~1.2.10].)}
 (1)
 Let ${\cal F}$ and ${\cal G}$ be
  an $\alpha$-twisted and a $\beta$ -twisted ${\cal O}_X$-module
  respectively,
  where $\alpha,\beta\in\check{C}^2_{\et}(X,{\cal O}_X^{\ast})$.
 Then
  ${\cal F}\otimes_{{\cal O}_X}{\cal G}$ is an $\alpha\beta$-twisted
   ${\cal O}_X$-module  and
  $\Homsheaf_{{\cal O}_X}({\cal F},{\cal G})$
   is an $\alpha^{-1}\beta$-twisted ${\cal O}_X$-module.
 In particular,
  if ${\cal F}$ and ${\cal G}$ are both
   $\alpha$-twisted ${\cal O}_X$-modules,
  then $\Homsheaf_{{\cal O}_X}({\cal F},{\cal G})$
   descends to an (ordinary/untwisted) ${\cal O}_X$-module,
   still denoted by $\Homsheaf_{{\cal O}_X}({\cal F},{\cal G})$,
   on $X$.

 (2)
 Let $f:X\rightarrow Y$ be a morphism of schemes$/{\Bbb C}$
  and $\alpha\in\check{C}^2_{\et}(Y,{\cal O}_Y^{\ast})$.
 Note that
  an $\alpha$-admissible \'{e}tale cover of $Y$ pulls back to
  an $f^{\ast}\alpha$-admissible \'{e}tale over of $X$ under $f$,
 through which the pull-back and push-forward of a related
  twisted sheaf can be defined.
 If ${\cal F}$ is an $\alpha$-twisted ${\cal O}_Y$-module on $Y$,
  then $f^{\ast}{\cal F}$ is an $f^{\ast}\alpha$-twisted
   ${\cal O}_X$-module on $X$.
 If ${\cal F}$ is an $f^{\ast}\alpha$-twisted
   ${\cal O}_X$-module on $X$,
  then $f_{\ast}{\cal F}$ is an $\alpha$-twisted
   ${\cal O}_Y$-module on $Y$.
\end{proposition}

\bigskip

\subsection{(General) Azumaya algebras over a scheme.}

\begin{definition}
{\bf [Azumaya algebra of rank $r$ over $X$].} {\rm
 An {\it Azumaya algebra ${\cal A}$ of rank\footnote{Note
                          that there are two conventions
                          in the literature:
                           rank as an ${\cal O}_X$-module vs.\
                           rank as an ${\cal O}_X$-algebra.
                          Here we take the latter convention.}
                          $r$}
  over a scheme $X/{\Bbb C}$
 is a locally free ${\cal O}_X$-algebra
 such that its fiber ${\cal A}\otimes_{{\cal O}_X} k(x)$
  at each closed point $x\in X$ is isomorphic to $\End({\Bbb C}^r)$
  ($=$ the $r\times r$ matrix algebra $M_r({\Bbb C})$ over ${\Bbb C}$)
  as ${\Bbb C}$-algebras.
}\end{definition}

\smallskip

\begin{proposition}
{\bf [local trivialization of Azumaya algebra].}
{\rm ([Mi: IV, Proposition 2.1].)}
 Let ${\cal A}$ be a sheaf of ${\cal O}_X$-algebras on $X$.
 The following statements are equivalent:
 \begin{itemize}
  \item[(1)]
   ${\cal A}$ is an Azumaya algebra of rank $r$ on $X$;

  \item[(2)]
   there is an \'{e}tale cover $U\rightarrow X$
    such that
    ${\cal A}\otimes_{{\cal O}_X}{\cal O}_U
     \simeq \Endsheaf_{{\cal O}_U}({\cal O}_U^{\oplus r})
     (=: M_r({\cal O}_U))$;

  \item[(3)]
   there is a flat cover $U\rightarrow X$
    such that
    ${\cal A}\otimes_{{\cal O}_X}{\cal O}_U
     \simeq \Endsheaf_{{\cal O}_U}({\cal O}_U^{\oplus r})$.
 \end{itemize}
\end{proposition}

\smallskip

\begin{definition-lemma}
{\bf [gerbe associated to Azumaya algebra].}
 {\it
  Given an Azumaya algebra ${\cal A}$ over $X$,
   the stack of trivializations of fibers of ${\cal A}$,
   defined by the assignment
   $$
    U\in \Scheme_X\;\longmapsto\;
     \mbox{\parbox[t]{28em}{the category with
      \begin{itemize}
       \item[$\cdot$]
        objects:\\{\rm
        the pairs $({\cal E}, a)$,
         where
          ${\cal E}$ is a locally free ${\cal O}_U$-module  and
          $a$ is an isomorphism
           $\Endsheaf_{{\cal O}_U}({\cal E})
            \stackrel{\sim}{\rightarrow}
             {\cal A}\otimes_{{\cal O}_X}{\cal O}_U$
           of ${\cal O}_U$-algebras,}

       \item[$\cdot$]
        morphism:\\{\rm
        a morphism $({\cal E}_1,a_1) \rightarrow ({\cal E}_2,a_2)$
         is an isomorphism
         $h:{\cal E}_1\stackrel{\sim}{\rightarrow}{\cal E}_2$
         such that the induced
          $h_{\ast}:\Endsheaf_{{\cal O}_U}({\cal E}_1)
                    \stackrel{\sim}{\rightarrow}
                    \Endsheaf_{{\cal O}_U}({\cal E}_2)$
          satisfies $a_2 \circ h_{\ast}=a_1$.}
      \end{itemize}
          } }
   $$
   is a gerbe, denoted by ${\cal X}_{\cal A}$,
   over $X$ with band ${\cal O}_X^{\ast}$.}
 {\rm
  We will call it the {\it gerbe over $X$ associated to ${\cal A}$}.
 The collection of objects $({\cal E}, a)$
  define a locally-free coherent
   ${\cal O}_{{\cal X}_{\cal A}}$-module ${\cal F}$
  on ${\cal X}_{\cal A}$.
 We will call it the {\it tautological fundamental module}
  on ${\cal X}_{\cal A}$.
 The pull-back of ${\cal A}$ on $X$ to ${\cal X}_{\cal A}$
  is canonically isomorphic to
  $\Endsheaf_{{\cal O}_{{\cal X}_{\cal A}}}({\cal F})$.
 We will call this the {\it tautological Azumaya algebra}
  over ${\cal X}_{\cal A}$.}
\end{definition-lemma}

\smallskip

Given an $\alpha$-twisted locally-free ${\cal O}_X$-module ${\cal E}$
 of rank $r$,
it follows from Proposition~1.1.7
 that ${\cal A}:=\Endsheaf_{{\cal O}_X}({\cal E})$
 is an Azumaya algebra over $X$.
Let $U\rightarrow X$ be an $\alpha$-admissible \'{e}tale cover of $X$.
Then a presentation of ${\cal E}$ on $U$ corresponds to a morphism
 $s:U\rightarrow {\cal X}_{\cal A}$.
In terms of this, ${\cal E}=s^{\ast}{\cal F}$,
 where ${\cal F}$
  is the tautological fundamental module on ${\cal X}_{\cal A}$.

\smallskip

\begin{definition-lemma}
{\bf [Brauer group $\Br(\,\bullet\,)$].} {\rm
 Two Azumaya algebras, ${\cal A}_1$ and ${\cal A}_2$,
  are said to be ({\it stably}) {\it equivalent}
  if there exist locally-free coherent ${\cal O}_X$-modules,
   ${\cal E}_1$ and ${\cal E}_2$,
  such that
   $\,{\cal A}_1\otimes_{{\cal O}_X}\Endsheaf_{{\cal O}_X}({\cal E}_1)\,
    \simeq\,
    {\cal A}_2\otimes_{{\cal O}_X}\Endsheaf_{{\cal O}_X}({\cal E}_2)\,$
   as ${\cal O}_X$-algebras.
 Denote the equivalence class of ${\cal A}$ by $[{\cal A}]$.
 Then, {\it
  the set of equivalence classes of Azumaya algebras over $X$
   form an abelian group, in notation $\Br(X)$,
   under
    $[{\cal A}_1]\cdot [{\cal A}_2]
     := [{\cal A}_1\otimes_{{\cal O}_X}{\cal A}_2]$,
    identity $= [{\cal O}_X]$,
    and $[{\cal A}]^{-1}:= [{\cal A}^{\circ}]$,
    where ${\cal A}^{\circ}$ is the opposite algebra\footnote{I.e.\
                      the Azumaya algebra
                      with the same ${\cal O}_X$-module ${\cal A}$
                      but with the reversed product $\cdot^{\prime}$
                      defined by $a_1\cdot^{\prime} a_2 := a_2\cdot a_1$
                      in the Azumaya algebra ${\cal A}$.}
    of ${\cal A}$}.
 $\Br(X)$ is called the {\it Brauer group} of $X$.
}\end{definition-lemma}

\smallskip

The set of isomorphism classes of Azumaya algebras of rank $r$ over $X$
 is given by the \'{e}tale cohomology group
 $H^1_{\et}(X,\PGL_r({\Bbb C}))$.
The exact sequence
 $$
  1\;\longrightarrow\; {\cal O}_X^{\ast}\;
     \longrightarrow\; \GL_r({\cal O}_X)\;
     \longrightarrow\; \PGL_r({\cal O}_X)\;
     \longrightarrow\; 1
 $$
 of sheaves on $X_{\et}$ defines an exact sequence of pointed-sets
 $$
  \cdots\; \longrightarrow\; H^1_{\et}(X,{\cal O}_X^{\ast})\;
  \longrightarrow\; H^1_{\et}(X,\GL_r({\cal O}_X))\;
  \longrightarrow\; H^1_{\et}(X,\PGL_r({\cal O}_X))\;
  \stackrel{d}{\longrightarrow}\; H^2_{\et}(X,{\cal O}_X^{\ast})\,.
 $$

\smallskip

\begin{theorem}
{\bf [$\Br(X)\subset H^2_{\et}(X,{\cal O}_X^{\ast})$].}
{\rm ([Mi: IV, Theorem 2.5].)}
 The connecting homomorphism
  $H^1_{\et}(X,\PGL_r({\cal O}_X))
   \stackrel{d}{\rightarrow} H^2_{\et}(X,{\cal O}_X^{\ast})$
  induces a canonical injective group-homomorphism
  $\Br(X)\hookrightarrow H^2_{\et}(X,{\cal O}_X^{\ast})$.
\end{theorem}

\bigskip

\section{Azumaya geometry and
         D-branes \`{a} la Polchinski-Grothendieck Ansatz revisited:
         the twist from a {\boldmath $B$}-field background.}

\subsection{Polchinski-Grothendieck Ansatz revisited
            with the \'{e}tale topology.}

\begin{flushleft}
{\bf Polchinski-Grothendieck Ansatz:
     Azumaya-type noncommutativity on D-branes.}
\end{flushleft}
Recall how the Polchinski-Grothendieck Ansatz for D-branes
 is reached in [L-Y1: Sec.~2.2].\footnote{Readers
               are referred to Polchinski [Pol2: vol.~I, Sec.~8.7]
                and [L-Y1: Sec.~2.2] for more thorough
                discussions and comparison.
               In this theme, we use as close notation to Polchinski
                as possible for a direct comparison.
               For all other parts of the work, we will use
                the more standard $X\rightarrow Y$
                to represent a D-brane (or D-brane world-volume) $X$
                that is mapped to a target-space(-time) $Y$.}
Consider a D-brane (or a D-brane world-volume) in a space(-time)
 that is geometrically realized
 as an embedded submanifold $f: Z\hookrightarrow M$
 in an open-string target space(-time) $M$.
The boundary of open-string world-sheets are mapped to $f(Z)$ in $M$.
Through this, open strings induce additional
 structures on $Z$, including
  a Chan-Paton bundle on $Z$ that supports the gauge field
    created from the vibrations of open-strings with end-points
    on $f(Z)$.
Let $\xi:=(\xi^a)_a$ be local coordinates on $Z$ and
  $X:=(X^a;X^{\mu})_{a,\mu}$ be local coordinates on $M$
  such that the embedding $f:Z\hookrightarrow M$ is locally
  expressed as
  $$
   X\; =\; X(\xi)\; =\; (X^a(\xi); X^{\mu}(\xi))_{a,\mu}\;
   =\; (\xi^a, X^{\mu}(\xi))_{a,\mu}\,;
  $$
 i.e., $X^a$'s (resp.\ $X^{\mu}$'s) are local coordinates along
       (resp.\ transverse to) $f(Z)$ in $M$.
This choice of local coordinates removes redundant degrees of freedom
 of the map $f$, and
$X^{\mu}=X^{\mu}(\xi)$ can be regarded as (scalar) fields on $Z$
 that collectively describes the postions/shapes/fluctuations
 of $Z$ in $M$ locally.
Here, both $\xi^a$'s, $X^a$'s, and $X^{\mu}$'s are ${\Bbb R}$-valued.
The gauge field on $Z$ is locally given by the connection
 $1$-form $A=\sum_a A_a(\xi)d\xi^a$ of a $U(1)$-bundle on $Z$.

When $r$-many such D-branes $Z$ are coincident, from the associated
 massless spectrum of (oriented) open strings with both end-points
 on $f(Z)$ one can draw the conclusion that
 \begin{itemize}
  \item[(1)]
   The gauge field $A=\sum_a A_a(\xi)d\xi^a$ on $Z$ is enhanced to
    $u(r)$-valued.

  \item[(2)]
   Each scalar field $X^{\mu}(\xi)$ on $Z$ is also enhanced
    to matrix-valued.
 \end{itemize}
Property (1) says that there is now a $U(r)$-bundle on $Z$.
To understand Property (2), one has two perspectives:
 \begin{itemize}
  \item[(A1)]
   [{\it coordinate tuple as point}]\hspace{1em}
   A tuple $(\xi^a)_a$ (resp.\ $(X^a; X^{\mu})_{a,\mu}$)
    represents a point on the world-volume $Z$ of the D-brane
    (resp.\ on the target space-time $M$).

  \item[(A2)]
   [{\it local coordinates as generating set of
         local functions}]\hspace{1em}
   Each local coordinate $\xi^a$ of $Z$ (resp.\ $X^a$, $X^{\mu}$ of $M$)
    is a local function on $Z$ (resp.\ on $M$)  and
   the local coordinates $\xi^a$'s
    (resp.\ $X^a$'s and $X^{\mu}$'s) together
    form a generating set of local functions on the world-volume $Z$
    of the D-brane (resp.\ on the target space-time $M$).
 \end{itemize}
While Aspect (A1) leads one to the anticipation of a noncommutative
  space from a noncommutatization of the target space-time $M$
  when probed by coincident D-branes,
 Aspect (A2) of Grothendieck leads one to a different/dual
   conclusion:
  a noncommutative space from a noncommutatization of
  the world-volume $Z$ of coincident D-branes,
 as follows.

Denote by ${\Bbb R}\langle \xi^a\rangle_{a}$
  (resp.\ ${\Bbb R}\langle X^a; X^{\mu}\rangle_{a, \mu}$)
 the local function ring on the associated local coordinate chart
 on $Z$ (resp.\ on $M$).
Then the embedding $f:Z\rightarrow M$,
  locally expressed as
  $X=X(\xi)=(X^a(\xi); X^{\mu}(\xi))_{a,\mu}=(\xi^a; X^{\mu}(\xi))$,
 is locally contravariantly equivalent to a ring-homomorphism
 $$
  f^{\sharp}\;:\;
   {\Bbb R}\langle X^a; X^{\mu}\rangle_{a, \mu}\;
   \longrightarrow\; {\Bbb R}\langle \xi^a\rangle_{a}\,,
  \hspace{1em}\mbox{generated by}\hspace{1em}
  X^a\;\longmapsto\; \xi^a\,,\;
  X^{\mu}\;\longmapsto\;X^{\mu}(\xi)\,.
 $$
When $r$-many such D-branes are coincident, $X^{\mu}(\xi)$'s become
 $M_r({\Bbb C})$-valued.
Thus, $f^{\sharp}$ is promoted to a new local ring-homomorphism:
 $$
  \hat{f}^{\sharp}\;:\;
   {\Bbb R}\langle X^a; X^{\mu}\rangle_{a, \mu}\;
   \longrightarrow\; M_r({\Bbb C}\langle \xi^a\rangle_{a})\,,
  \hspace{1em}\mbox{generated by}\hspace{1em}
  X^a\;\longmapsto\; \xi^a\cdot{\mathbf 1}\,,\;
  X^{\mu}\;\longmapsto\;X^{\mu}(\xi)\,.
 $$
Under Grothendieck's contravariant local equivalence of function rings
 and spaces, $\hat{f}^{\sharp}$ is equivalent to saying that we have
 now a map $\hat{f}: Z_{\scriptsizenoncommutative}\rightarrow M$.
Thus, the D-brane-related noncommutativity in Polchinski's treatise
 re-read from the viewpoint of
 Grothendieck implies the following ansatz:

\bigskip

\noindent
{\bf Polchinski-Grothendieck Ansatz [D-brane: noncommutativity].}
{\it
 The world-volume of a D-brane carries a noncommutative structure
 locally associated to a function ring of the form $M_r(R)$,
 where $r\in {\Bbb Z}_{\ge 1}$ and
   $M_r(R)$ is the $r\times r$ matrix ring over $R$.
} 

\bigskip

\noindent
Note that $R$ can be either commutative or noncommutative,
cf.\ Remark~5.1.9.

\bigskip

Cf.~[L-L-S-Y: {\sc Figure}~1-2].

%
%
%

\bigskip

\bigskip

\begin{flushleft}
{\bf Polchinski-Grothendieck Ansatz
     with the \'{e}tale topology adaptation.}
\end{flushleft}
In the {\it smooth differential-geometric setting} of Polchinski,
 the word ``{\it locally}" in the ansatz means
 ``locally in the {\it $C^{\infty}$-topology}".
This can be generalized to adapt the ansatz to fit various settings:
 ``locally" in the {\it analytic} (resp.\ {\it Zariski}) {\it topology}
 for the {\it holomorphic} (resp.\ {\it algebro-geometric}) {\it setting}.
These are enough to study D-branes in a space(-time) without
 a background $B$-field.
The Azumaya structure sheaf ${\cal O}_Z^{A\!z}$
 that encodes the matrix-type noncommutative structure on $Z$
 in these cases is then of the form $\Endsheaf_{{\cal O}_Z}({\cal E})$
 with ${\cal E}$ the Chan-Paton module,
 a locally free ${\cal O}_Z$-module of rank $r$
 on which the Azumaya ${\cal O}_Z$-algebra ${\cal O}_Z^{A\!z}$
 acts tautologically
 as a simple/fundamental (left) ${\cal O}_Z^{A\!z}$-module.
This is the case studied in the previous part
 [L-Y1], [L-L-S-Y], [L-Y2], [L-Y3] of the project.
${\cal O}_Z^{A\!z}$ in this case corresponds to the zero-class
 in the Brauer group $\Br(Z)$ of $Z$.

On pure mathematical ground, one can further adapt the ansatz
 for $Z$ equipped with any Grothendieck topolgy/site.
On string-theoretic ground, as recalled in Sec.~0,
when a background $B$-field $B$ on $M$ is turned on,
 the Chan-Paton module ${\cal E}$ on $Z$ becomes twisted and
 is no longer an honest sheaf of ${\cal O}_Z$-modules on $Z$.
The interpretation of ``{\it locally}" in the ansatz
  in the sense of (small) {\it \'{e}tale topology} on $Z$
 becomes forced upon us.
This corresponds to the case
 when the Azumaya structure sheaf ${\cal O}_Z^{A\!z}$ on $Z$
 represents a non-zero class in $\Br(Z)$.

We now turn to the algebro-geometric aspect of D-branes,
 following {\it Polchinski-Grothendieck Ansatz}
 but {\it with} this {\it \'{e}tale topology adaptation}
  on the D-brane or D-brane world-volume.

\bigskip

\subsection{D-branes in a {\boldmath $B$}-field background
   as morphisms from Azumaya schemes with a twisted fundamental module.}

Recall the twisting effect of $B$-field from string theory
 highlighted in Introduction.
We now study D-branes in a $B$-field background
 along the line of the adapted Polchinski-Grothendieck Ansatz.
Except the additional involvement of \'{e}tale topology,
  twisted sheaves, and the matching of twists,
 the setting/study in [L-Y1: Sec.~1] and [L-L-S-Y: Sec.~2]
  carries over directly to the current situation.

\bigskip

\begin{flushleft}
{\bf D-branes in a {\boldmath $B$}-field background.}
\end{flushleft}
\begin{definition}
{\bf [Azumaya scheme with a fundamental module].}
{\rm
 An {\it Azumaya scheme with a fundamental module in class $\alpha$}
  is a tuple
  $$
   (X^{A\!z},{\cal E})\;
    :=\; (X,\,
          {\cal O}_X^{A\!z}
            = \Endsheaf_{{\cal O}_X}({\cal E}),\,
          {\cal E})\,,
  $$
  where
   $X=(X,{\cal O}_X)$ is a (Noetherian) scheme (over ${\Bbb C}$),
   $\alpha\in \check{C}^2_{\et}(X,{\cal O}_X^{\ast})$ represents
    a class $[\alpha]\in\Br(X)\subset H^2_{\et}(X,{\cal O}_X^{\ast})$,
    and
   ${\cal E}$ is a locally-free coherent $\alpha$-twisted
    ${\cal O}_X$-module on $X$.
 A {\it commutative surrogate} of $(X^{A\!z},{\cal E})$
  is a scheme $X_{\cal A}:=\boldSpec{\cal A}$,
  where
   ${\cal O}_X\subset {\cal A}\subset \Endsheaf_{{\cal O}_X}({\cal E})$
   is an inclusion sequence of commutative ${\cal O}_X$-subalgebras
   of $\Endsheaf_{{\cal O}_X}({\cal E})$.
 Let $\pi: X_{\cal A}\rightarrow X$
  be the built-in dominant finite morphism.
 Then
  $\cal E$ is tautologically a $\pi^{\ast}\alpha$-twisted
   ${\cal O}_{X_{\cal A}}$-module on $X_{\cal A}$,
   denoted by $_{{\cal O}_{X_{\cal A}}}{\cal E}$.
 We say that $X^{A\!z}$ is an {\it Azumaya scheme of rank $r$}
   if ${\cal E}$ has rank $r$ and
  that it is a {\it nontrivial} (resp.\ {\it trivial})
   Azumaya scheme if $[\alpha]\ne 0$ (resp.\ $[\alpha]=0$).
}\end{definition}

\smallskip

\begin{remark}
{\it $[\,$the twisted sheaf $_{{\cal O}_{X_{\cal A}}}{\cal E}$
          on $X_{\cal A}$$\,]$.}
{\rm
 Explicitly, let
  $p:U\rightarrow X$
   be an $\alpha$-admissible \'{e}tale cover of $X$
   on which the $\alpha$-twisted ${\cal O}_X$-module ${\cal E}$
   is represented as an ordinary ${\cal O}_U$-module  and
  $p_{\cal A}: U_{\cal A}:= U\times_X X_{\cal A}\rightarrow X_{\cal A}$
   be the pull-back \'{e}tale cover of $X_{\cal A}\,$:
   $$
    \xymatrix{
     & U_{\cal A} \ar[d]_{\widetilde{\pi}} \ar[rr]^-{p_{\cal A}}
                  && X_{\cal A} \ar[d]^{\pi}\\
     & U \ar[rr]^-{p}           && X     &.
    }
   $$
 Then
  ${\cal O}_{U}\subset p^{\ast}{\cal A}  \subset
    p^{\ast}{\cal O}_X^{A\!z}=\Endsheaf_{{\cal O}_U}({\cal E})
    = {\cal O}_U^{A\!z}$
   is a sequence of ${\cal O}_U$-subalgebra inclusions
    and
  $\widetilde{\pi}:U_{\cal A}\rightarrow U$ is nothing but
   the commutative surrogate $\boldSpec(p^{\ast}{\cal A})\rightarrow U$
   of $U^{A\!z}=(U,{\cal O}_U^{A\!z})$.
  In terms of this, ${\cal E}$ is canonically
   a $p^{\ast}{\cal A}={\cal O}_{U_{\cal A}}$-module,
   which defines then
    the $\pi^{\ast}\alpha$-twisted ${\cal O}_{X_{\cal A}}$-module
    $_{{\cal O}_{X_{\cal A}}}{\cal E}$ on $X_{\cal A}$
    in the above definition.
}\end{remark}

\smallskip

\begin{remark}
{\it $[\,$noncommutative-space viewpoint$\,]$.} {\rm
 It is instructive to think of
  ${\cal O}_X^{A\!z}$ as the structure sheaf ${\cal O}$ of
    a noncommutative space $\Space{\cal O}_X^{A\!z}$,
  $\Space{\cal O}_U^{A\!z}$ as an \'{e}tale cover of
   $\Space{\cal O}_X^{A\!z}$,
  ${\cal E}$ as an ordinary sheaf of (left) modules
   on $\Space{\cal O}_U^{A\!z}$
   that defines a (left) twisted ${\cal O}$-module
   on $\Space{\cal O}_X^{A\!z}$,
   and
  that there is a dominant morphism
   $\Space{\cal O}_X^{A\!z}\rightarrow X_{\cal A}$,
   under which ${\cal E}$ on $\Space{\cal O}_X^{A\!z}$
    is pushed forward to the $\pi^{\ast}\alpha$-twisted
    ${\cal O}_{X_{\cal A}}$-module
    $_{{\cal O}_{X_{\cal A}}}{\cal E}$ on $X_{\cal A}$.
}\end{remark}

\smallskip

Let
 $Y$ be a (commutative, Noetherian) scheme/${\Bbb C}$  and
 $\alpha_B\in \check{C}^2_{\et}(Y,{\cal O}_Y^{\ast})$
  be the \'{e}tale \v{C}ech cocycle associated to a fixed $B$-field
  on $Y$.
Then a proto-typical definition of D-branes (of B/holomorphic type)
 in $(Y,B)$ is given by morphisms
 $\varphi:(X^{A\!z},{\cal E})\rightarrow (Y,\alpha_B)$,
 defined as follows.

\smallskip

\begin{definition}
{\bf [morphism from Azumaya scheme with fundamental module to
      {\boldmath $B$}-field background].}
{\rm
 Let $(X^{A\!z},{\cal E})$ be an Azumaya scheme with a fundamental module
  in the class $\alpha\in\check{C}^2_{\et}(X,{\cal O}_X^{\ast})$.
 Then,
 a {\it morphism} from $(X^{A\!z},{\cal E})$ to $(Y,\alpha_B)$,
  in notation $\varphi: (X^{A\!z},{\cal E})\rightarrow (Y,\alpha_B)$,
  is a pair
  $$
   ({\cal O}_X \subset {\cal A}_{\varphi}
               \subset {\cal O}_X^{\Azscriptsize}\;,\;
     f_{\varphi}: X_{\varphi}:=\boldSpec{\cal A}_{\varphi}
                  \rightarrow Y)\,,
  $$
  where
  \begin{itemize}
   \item[$\cdot$]
    ${\cal A}_{\varphi}$ is a commutative ${\cal O}_X$-subalgebra
    of ${\cal O}_X^{A\!z}$,

   \item[$\cdot$]
    $f_{\varphi}:X_{\varphi} \rightarrow Y$
    is a morphism of (commutative) schemes,
  \end{itemize}
  that satisfies the following properties:
  \begin{itemize}
   \item[(1)]
   ({\it minimal property of $X_{\varphi}$})$\;$
    there exists no ${\cal O}_X$-subalgebra
     ${\cal O}_X \subset {\cal A}^{\prime}\subset {\cal A}_{\varphi}$
    such that
     $f_{\varphi}$ factors as the composition of morphisms
     $X_{\varphi}
      \rightarrow \boldSpec{\cal A}^{\prime} \rightarrow Y$;

   \item[(2)]
   ({\it matching of twists on $X_{\varphi}$})$\;$
    let $\pi_{\varphi}:X_{\varphi}\rightarrow X$
     be the built-in finite dominant morphism,
    then $\pi_{\varphi}^{\ast}\alpha=f_{\varphi}^{\ast}\alpha_B$
     in $\check{C}^2_{\et}(X_{\varphi},{\cal O}_{X_{\varphi}}^{\ast})$.
  \end{itemize}
 $X_{\varphi}$ is called the
  {\it surrogate of $X^{\Azscriptsize}$ associated to $\varphi$}.
 Condition (2) implies that
  $\varphi_{\ast}{\cal E}
   := f_{\varphi\,\ast}(_{{\cal O}_{X_{\varphi}}}{\cal E})$
  is an $\alpha_B$-twisted ${\cal O}_Y$-module on $Y$,
  supported on
  $\Image(\varphi) :=\varphi(X^{A\!z})
   := f_{\varphi}(X_{\varphi})$,
   where the last is the usual scheme-theoretic image
    of $X_{\varphi}$ under $f_{\varphi}$.\footnote{In
                              other words,
                              a morphism from $(X^{A\!z},{\cal E})$
                               to $(Y,\alpha_B)$ is a usual morphism
                               $\varphi: X^{A\!z} \rightarrow Y$
                               from the (possibly nontrivial)
                               Azumaya scheme $X^{A\!z}$ to $Y$
                              subject to the twist-matching Condition (2)
                              so that $\varphi_{\ast}{\cal E}$
                               remains a twisted sheaf
                               in a way that is compatible with
                               the $B$-field background on $Y$.}

 Given two morphisms
    $\varphi_1:(X_1^{A\!z},{\cal E}_1)\rightarrow (Y,\alpha_B)$ and
    $\varphi_2:(X_2^{A\!z},{\cal E}_2)\rightarrow (Y,\alpha_B)$,
 a {\it morphism} $\varphi_1\rightarrow\varphi_2$
  from $\varphi_1$ to $\varphi_2$ is a pair $(h, \widetilde{h})$,
  where
  \begin{itemize}
   \item[$\cdot$]
    $h:X_1\rightarrow X_2$ is an isomorphism of schemes
      with $h^{\ast}\alpha_2=\alpha_1$,
      where $\alpha_i$ is the underlying class of ${\cal E}_i$
       in $\check{C}^2_{\et}(X_i,{\cal O}_{X_i}^{\ast})$;

   \item[$\cdot$]
    $\widetilde{h}:{\cal E}_1
                    \stackrel{\sim}{\rightarrow} h^{\ast}{\cal E}_2$
     be an isomorphism of twisted sheaves on $X_1$
     that satisfies
     \begin{itemize}
      \item[$\cdot$]
       $\widetilde{h}:
        {\cal A}_{\varphi_1}
        \stackrel{\sim}{\rightarrow} h^{\ast}{\cal A}_{\varphi_2}$,

      \item[$\cdot$]
       the following diagram commutes
       \begin{eqnarray*}
        \xymatrix{
          X_{\varphi_2}\ar[drr]^{f_{\varphi_2}}\ar[d]_{\widehat{h}}
                                                           &&      \\
          X_{\varphi_1}\ar[rr]^{f_{\varphi_1}}             && Y\; .\\
        }
       \end{eqnarray*}
     \end{itemize}
    Here,
    we denote both of
     the induced isomorphisms,
      $\, {\cal O}_{X_1}^{A\!z} \stackrel{\sim}{\rightarrow}
                              h^{\ast}{\cal O}_{X_2}^{A\!z}$
       and
      ${\cal A}_{\varphi_1}
        \stackrel{\sim}{\rightarrow} h^{\ast}{\cal A}_{\varphi_2}\,$,
      of ${\cal O}_{X_1}$-algebras still by $\widetilde{h}$  and
    $\widehat{h}:X_{\varphi_2}
                 \stackrel{\sim}{\rightarrow} X_{\varphi_1}$
     is the scheme-isomorphism associated to
    $\widetilde{h}: {\cal A}_{\varphi_1}
      \stackrel{\sim}{\rightarrow} h^{\ast}{\cal A}_{\varphi_2}$.
  \end{itemize}

 This defines the category $\MorphismCategory_{A\!z^f}(Y,\alpha_B)$
  of morphisms
  from Azumaya schemes with a fundamental module
  to $(Y,\alpha_B)$.
}\end{definition}

\smallskip

\begin{definition}
{\bf [D-brane and Chan-Paton module].} {\rm
 Following the previous Definition,
  $\varphi(X^{A\!z})$
   is called the {\it image D-brane} on $(Y,\alpha_B)$  and
  $\varphi_{\ast}{\cal E}$ the {\it Chan-Paton module}
   on the image D-brane.
 Similarly, for {\it image D-brane world-volume}
  if $X$ is served as a (Wicked-rotated) D-brane world-volume.
}\end{definition}

\smallskip

\begin{remark}
{\it $[\,$fundamental vs.\ solitonic D-brane$\,]$.} {\rm
 The setting here treats D-branes in string theory
  more as a fundamental/soft extended object.
 For {\it solitonic}/hard D-brane in space-time,
  one may require in addition that
  $f_{\varphi}:X_{\varphi}\rightarrow Y$ be an embedding.
}\end{remark}

\bigskip

\begin{flushleft}
{\bf Azumaya without Azumaya and morphisms without morphisms.}
\end{flushleft}
Similar to the case of trivial Azumaya curves
 studied in [L-L-S-Y: Sec.~2],
Definition~2.2.4
 has an equivalent version in terms of twisted sheaves on
 the related product space as follows.

Given a morphism
  $\varphi:(X^{A\!z}, {\cal E})\rightarrow (Y,\alpha_B)$
  as in Definition~2.2.4,
 the minimal property of the surrogate $X_{\varphi}$ of $X^{A\!z}$
  associated to $\varphi$ implies that
 $(\pi_{\varphi},f_{\varphi}):X_{\varphi}\rightarrow X\times Y$
  embeds $X_{\varphi}$ in $X\times Y$
  as a subscheme $\Gamma_{\varphi}$.
Let $\pr_1:X\times Y\rightarrow X$ and $\pr_2:X\times Y\rightarrow Y$
 be the projection maps.
Then
 the $\pi_{\varphi}^{\ast}{\alpha}$-twisted
  ${\cal O}_{X_{\varphi}}$-module
  $_{{\cal O}_{X_{\varphi}}}{\cal E}$ on $X_{\varphi}$
 is pushed forward to a $\pr_1^{\ast}\alpha$-twisted
  ${\cal O}_{X\times Y}$-module $\widetilde{\cal E}$ on $X\times Y$
  that is supported on $\Gamma_{\varphi}$.
The matching condition of twists on $\varphi$ says that
 $\pr_1^{\ast}\alpha=\pr_2^{\ast}\alpha_B$ on $\Gamma_{\varphi}$.
By construction, $\widetilde{\cal E}$ on $X\times Y$
 is flat over $X$ with relative length $r$.

Conversely, given
 a $(\alpha,\alpha_B)
        \in \check{C}^2_{\et}(X,{\cal O}_X^{\ast})
                     \times \check{C}^2_{\et}(Y,{\cal O}_Y^{\ast})$
  as before  and
 a coherent $\pr_1^{\ast}\alpha$-twisted ${\cal O}_{X\times Y}$-module
  $\widetilde{\cal E}$ on $X\times Y$
  that satisfies the following two conditions:
  \begin{itemize}
   \item[(a)]
    $\widetilde{\cal E}$ is flat over $X$ with relative length $r$;

   \item[(b)]
    $\pr_1^{\ast}\alpha=\pr_2^{\ast}\alpha_B$
    on $\Supp\widetilde{\cal E} =:\Gamma$.
  \end{itemize}
Then,
 ${\cal E} := \pr_{1\ast}\widetilde{\cal E}$
  is an $\alpha$-twisted ${\cal O}_X$-module on $X$.
 This defines an Azumaya scheme with a fundamental module, i.e.\
  $(X,{\cal O}_X^{A\!z}:=\Endsheaf_{{\cal O}_X}({\cal E}),{\cal E})$,
  over $X$ in the class $\alpha$.
The defining ${\cal O}_{\Gamma}$-algebra homomorphism
 ${\cal O}_{\Gamma}
  \rightarrow \Endsheaf_{{\cal O}_{\Gamma}}(\widetilde{\cal E})$
 realizes ${\cal O}_{\Gamma}$ as an ${\cal O}_X$-algebra ${\cal A}$
 that fits into ${\cal O}_X\subset {\cal A}\subset {\cal O}_X^{A\!z}$
  canonically.
By construction,
 $X_{\cal A}:=\boldSpec{\cal A} \simeq \Gamma$ canonically  and
 the restriction
   $\xymatrix{X & \Gamma \ar[l]_-{\prscriptsize_1}
                       \ar[r]^-{\prscriptsize_2}  & Y}$
   of the projection maps
  defines morphisms
   $\xymatrix{X & X_{\cal A} \ar[l]_-{\pi} \ar[r]^-{f}  & Y}$
   that satisfies both Condition (1) (minimal property) and
   Condition (2) (matching of twists)
   in Definition~2.2.4.
$f$ thus defines a morphism
 $\varphi:(X^{A\!z},{\cal E})\rightarrow (Y,\alpha_B)$
 with $X_{\varphi}=X_{\cal A}$ and $f_{\varphi}=f$.

Let $\CohCategory_{\productscriptsize^0}(Y,\alpha_B)$
 be the category with objects
 coherent twisted modules $\widetilde{\cal E}$ on a product $X\times Y$
  that is flat over $X$ with relative dimension $0$ and
   satisfies Condition (b) above.
A {\it morphism}
  $\widetilde{\cal E}_1\rightarrow \widetilde{\cal E}_2$
  from $\widetilde{\cal E}_1$ on $X_1\times Y$
  to $\widetilde{\cal E}_2$ on $X_2\times Y$
 is a pair $(h,\widetilde{h})$
 where
  \begin{itemize}
   \item[$\cdot$]
    $h:X_1\rightarrow X_2$ is an isomorphism of schemes
      with $h^{\ast}\alpha_2=\alpha_1$,
      where $\alpha_i$ is the underlying class
       in $\check{C}^2_{\et}(X_i,{\cal O}_{X_i}^{\ast})$
       in question;

   \item[$\cdot$]
    denote the induced isomorphism
     $X_1\times Y\stackrel{\sim}{\rightarrow} X_2\times Y$
     also by $h$,
    then
    $\widetilde{h}:\widetilde{\cal E}_1
      \stackrel{\sim}{\rightarrow} h^{\ast}\widetilde{\cal E}_2$
     is an isomorphism of $\alpha_1$-twisted ${\cal O}_{X_1}$-modules
     on $X_1$.
  \end{itemize}
The discussion above defines two functors
 $$
  \xymatrix{
   \MorphismCategory_{A\!z^f}(Y,\alpha_B) \ar@<.4ex>[rr]^-{F}
    && \CohCategory_{\productscriptsize^0}(Y,\alpha_B)
        \ar@<.4ex>[ll]^-{G}\,.
  }
 $$

\smallskip

\begin{lemma}
{\bf [Azumaya without Azumaya, morphisms without morphisms].}
 $(F,G)$ defines an equivalence of categories
 $\MorphismCategory_{A\!z^f}(Y,\alpha_B)$ and
 $\CohCategory_{\productscriptsize^0}(Y,\alpha_B)$.
\end{lemma}

\bigskip

\begin{flushleft}
{\bf The description in terms of morphisms from Azumaya gerbes
     with a fundamental module to a target gerbe.\footnote{This
                              theme is written with
                             the work [Sh3] of Eric Sharpe,
                               particularly
                               [Sh3: Sec.~6.3 {\it D-brane ``bundles"}]
                                and
                               [Sh3: Sec.~7 {\it Conclusions}]
                                concerning {\it the equivalence of
                                ``turning on the $B$-field" and
                                ``compactifying on a generalized space
                                (i.e.\ a gerbe or a sheaf on a gerbe)"},
                              also in mind  and, hence,
                             goes with a hidden subtitle:
                              {\it Sharpe vs.\ Polchinski-Grothendieck}.
                             Readers are highly recommended to read
                              ibidem alongside.
                             We thank him for comments on [L-Y1]
                              and sharing with us his insights
                              on various subtle issues in string theory
                              in fall 2007.}}
\end{flushleft}
Given an Azumaya algebra
  ${\cal A}=\Endsheaf_{{\cal O}_X}({\cal E})$ over $X$,
 where ${\cal E}$
  is an $\alpha$-twisted locally free ${\cal O}_X$-module
   defined on an $\alpha$-admissible \'{e}tale cover
   $p:U\rightarrow X$ of $X$,
recall
 the ${\cal O}_X^{\ast}$-gerbe ${\cal X}_{\cal A}$ over $X$  and
 the tautological fundamental module ${\cal F}$ on ${\cal X}_{\cal A}$.
Then ${\cal E}$ defines an atlas
 $\breve{p}:U\rightarrow {\cal X}_{\cal A}$
  of the algebraic stack ${\cal X}_{\cal A}$,
 with $\breve{p}^{\ast}{\cal F}={\cal E}$.
The pull-back $\breve{\cal A}$ of ${\cal A}$ to ${\cal X}_{\cal A}$
 is canonically isomorphic to
 ${\cal O}^{A\!z}_{{\cal X}_{\cal A}}
         :=\Endsheaf_{{\cal O}_{{\cal X}_{\cal A}}}({\cal F})$
  on ${\cal X}_{\cal A}$.
The latter defines the tautological trivial Azumaya structure
 on ${\cal X}_{\cal A}$ with the tautological fundamental module
 ${\cal F}$.
\begin{itemize}
 \item[$\cdot$] [{\it notation}]\hspace{1em}
  We will denote
   the gerbe ${\cal X}_{\cal A}$ by ${\cal X}$  and
   the Azumaya algebraic stack with a fundamental module
    $({\cal X}_{\cal A}, {\cal O}_{{\cal X}_{\cal A}}^{A\!z}, {\cal F})$
    by $({\cal X}^{A\!z},{\cal F})$
  for simplicity in the following discussion.
\end{itemize}
The notion of {\it surrogates} of an Azumaya scheme generalizes
 directly to that of an Azumaya algebraic stack   and
the notion of {\it morphisms} from a trivial Azumaya scheme
 with a fundamental module generalizes directly -
   after the combination with the notion of morphisms of algebraic stacks -
  to that from a trivial Azumaya algebraic stack with a fundamental module.

In our case,
let $\varphi:(X^{A\!z},{\cal E})\rightarrow (Y,\alpha_B)$
 be a morphism, specified by a pair
 $$
  ({\cal O}_X \subset {\cal A}_{\varphi}  \subset
     {\cal O}_X^{\Azscriptsize}:=\Endsheaf_{{\cal O}_X}({\cal E})\;,\;
   f_{\varphi}: X_{\varphi}:=\boldSpec{\cal A}_{\varphi} \rightarrow Y)\,.
 $$
Let ${\cal X}_{\varphi}= {\cal X}\times_X X_{\varphi}$.
This is a gerbe over $X_{\varphi}$
  in the class $\pi_{\varphi}^{\ast}\alpha$,
 with an atlas
  $$
   \breve{p}_{\varphi}\; :\;
     U_{\varphi}:= U_{{\cal A}_{\varphi}} := U\times_X X_{\varphi}\;
                                   \longrightarrow\; {\cal X}_{\varphi}
  $$
  induced by
  $p_{\varphi}:=\pi_{\varphi}^{\ast}p:U_{\varphi}\rightarrow X_{\varphi}$.
The ${\cal O}_X$-subalgebra ${\cal A}_{\varphi}$ of ${\cal O}_X^{A\!z}$
 induces an ${\cal O}_{\cal X}$-subalgebra $\breve{\cal A}_{\varphi}$
  of ${\cal O}_{\cal X}^{A\!z}$.
This defines a surrogate ${\cal X}_{\breve{A}_{\varphi}}$
 of ${\cal X}^{A\!z}$ that is precisely ${\cal X}_{\varphi}$.
Let
 ${\cal Y}={\cal Y}_{\alpha_B}$ be an ${\cal O}_Y$-gerbe over $Y$
  that represents $\alpha_B$  and
 $\breve{q}: V\rightarrow {\cal Y}$ be an atlas of ${\cal Y}$
  with the following properties:
  \begin{itemize}
   \item[(1)]
    the underlying $q: V\rightarrow Y$
    gives an $\alpha_B$-admissible \'{e}tale cover of $Y$;

   \item[(2)]
    $\underline{\Isom}(\breve{q},\breve{q}) := V\times_{\cal Y}V$
    has a global section over $V\times_Y V$.
  \end{itemize}
(Note that
 ${\cal Y}$ is non-empty on each connected component of $V$.)
Let
 $p_{\varphi}^{\prime}:U_{\varphi}^{\prime}\rightarrow X_{\varphi}$
  be an \'{e}tale refinement of $p_{\varphi}$ of $X_{\varphi}$
 so that:
 \begin{itemize}
  \item[(1)]
   $p_{\varphi}^{\prime}$ refines also the \'{e}tale cover
    $f_{\varphi}^{\ast}q:X_{\varphi}\times_Y V\rightarrow X_{\varphi}$
    of $X_{\varphi}$;

  \item[(2)]
   $\underline{\Isom}
    (\breve{p}_{\varphi}^{\prime},\breve{p}_{\varphi}^{\prime})  :=
    U_{\varphi}^{\prime}\times_{{\cal X}_{\varphi}} U_{\varphi}^{\prime}$
   has a global section over
    $U_{\varphi}^{\prime}\times_{X_{\varphi}} U_{\varphi}^{\prime}$,
   where
    $\breve{p}_{\varphi}^{\prime}:
             U_{\varphi}^{\prime} \rightarrow {\cal X}_{\varphi}$
    is the new atlas of ${\cal X}_{\varphi}$
    associated to the refinement
     $U_{\varphi}^{\prime}\rightarrow U_{\varphi}$.
 \end{itemize}
Then, one has the following diagram
 $$
  \xymatrix{
     **[r]U_{\varphi}^{\prime}\ar[d]_{p_{\varphi}^{\prime}}
          \ar[rr]^-{\hat{f}_{\varphi}}        &&  V\ar[d]^q \\
     **[r]X_{\varphi} \ar[rr]^-{f_{\varphi}} &&  **[r]Y\,,
  }
 $$
 where $\hat{f}_{\varphi}$ is the composition
  $U_{\varphi}^{\prime}\rightarrow X_{\varphi}\times_YV\rightarrow V$.

Fix
 a global section of
  the ${\cal O}_{X_{\varphi}}^{\ast}$-torsor
   $U_{\varphi}^{\prime}\times_{{\cal X}_{\varphi}} U_{\varphi}^{\prime}$
   over $U_{\varphi}^{\prime}\times_{X_{\varphi}}U_{\varphi}^{\prime}$
     and
 a global section of
  the ${\cal O}_Y^{\ast}$-torsor
   $V\times_{\cal Y}V$ over $V\times_Y V$.
This trivializes
 the ${\cal O}_{\bullet}^{\ast}$-torsors
  $(U_{\varphi}^{\prime}
                \times_{{\cal X}_{\varphi}} U_{\varphi}^{\prime})/
         (U_{\varphi}^{\prime}\times_{X_{\varphi}}U_{\varphi}^{\prime})$
     and
  $(V\times_{\cal Y}V)/(V\times_Y V)$,
 the ${\cal O}_{\bullet}^{\ast}\times{\cal O}_{\ast}^{\ast}$-torsors
  $(U_{\varphi}^{\prime}
     \times_{{\cal X}_{\varphi}} U_{\varphi}^{\prime}
     \times_{{\cal X}_{\varphi}} U_{\varphi}^{\prime})/
         (U_{\varphi}^{\prime}
              \times_{X_{\varphi}} U_{\varphi}^{\prime}
              \times_{X_{\varphi}} U_{\varphi}^{\prime})$
     and
  $(V\times_{\cal Y}V\times_{\cal Y}V)/(V\times_Y V\times_Y V)$,
  $\,\cdots\,,\,$ etc.\
and
it follows that $\hat{f}_{\varphi}$ lifts to
 a commutative diagram of multi-arrows:\footnote{Here,
                        it is understood that
                         a commutative diagram applies only to
                         a square in the tower
                         with same-type projection maps
                         for its vertical arrows.
                        E.g.\
                         $\hat{f}_{\varphi}^{\,(1)}
                           \circ p_{\varphi, 13}^{\prime}
                           = q_{13}\circ \hat{f}_{\varphi}^{\,(2)}$.}
 $$
  \xymatrix{
   \hspace{1em}\vdots\hspace{1em}
    \ar[rr]^-{\hat{f}_{\varphi}^{\,(\bullet)}}
    \ar@<-1.2ex>[d] \ar@<-.4ex>[d] \ar@<.4ex>[d] \ar@<1.2ex>[d]
     && \hspace{1em}\vdots\hspace{1em}
        \ar@<-1.2ex>[d]\ar@<-.4ex>[d]\ar@<.4ex>[d]\ar@<1.2ex>[d]\\
   U_{\varphi}^{\prime}
     \times_{{\cal X}_{\varphi}} U_{\varphi}^{\prime}
     \times_{{\cal X}_{\varphi}} U_{\varphi}^{\prime}
       \ar[rr]^-{\hat{f}_{\varphi}^{\,(2)}}
       \ar@<-.8ex>[d] \ar[d] \ar@<.8ex>[d]
     && V\times_{\cal Y}V\times_{\cal Y}V
       \ar@<-.8ex>[d] \ar[d] \ar@<.8ex>[d] \\
   U_{\varphi}^{\prime}\times_{{\cal X}_{\varphi}} U_{\varphi}^{\prime}
    \ar[rr]^-{\hat{f}_{\varphi}^{\,(1)}}
    \ar@<-.4ex>[d] \ar@<.4ex>[d]
     && V\times_{\cal Y}V \ar@<-.4ex>[d] \ar@<.4ex>[d] \\
   **[r]U_{\varphi}^{\prime}\ar[d]_{p_{\varphi}^{\prime}}
        \ar[rr]^-{\hat{f}_{\varphi}}        &&  V\ar[d]^q \\
   **[r]X_{\varphi} \ar[rr]^-{f_{\varphi}} &&  Y
  }
 $$
 that covers - indeed trivialized trivial torsor over -
 the $\hat{f}_{\varphi}$-induced tautological tower
 $$
  \xymatrix{
   \hspace{1em}\vdots\hspace{1em}
    \ar[rr] \ar@<-1.2ex>[d] \ar@<-.4ex>[d] \ar@<.4ex>[d] \ar@<1.2ex>[d]
     && \hspace{1em}\vdots\hspace{1em}
        \ar@<-1.2ex>[d]\ar@<-.4ex>[d]\ar@<.4ex>[d]\ar@<1.2ex>[d]\\
   U_{\varphi}^{\prime}
     \times_{X_{\varphi}} U_{\varphi}^{\prime}
     \times_{X_{\varphi}} U_{\varphi}^{\prime}
      \ar[rr] \ar@<-.8ex>[d] \ar[d] \ar@<.8ex>[d]
     && V\times_YV\times_Y V \ar@<-.8ex>[d] \ar[d] \ar@<.8ex>[d] \\
   U_{\varphi}^{\prime}\times_{X_{\varphi}} U_{\varphi}^{\prime}
    \ar[rr] \ar@<-.4ex>[d] \ar@<.4ex>[d]
     && V\times_Y V \ar@<-.4ex>[d] \ar@<.4ex>[d] \\
   **[r]U_{\varphi}^{\prime}\ar[d]_{p_{\varphi}^{\prime}}
        \ar[rr]^-{\hat{f}_{\varphi}}        &&  V\ar[d]^q \\
   **[r]X_{\varphi} \ar[rr]^-{f_{\varphi}} &&  **[r]Y\;.
  }
 $$
A standard decent-data argument implies then that:

\smallskip

\begin{lemma}
{\bf [presentation as morphism from Azumaya gerbe].}
 $f_{\varphi}:X_{\varphi}\rightarrow Y$
  induces a morphism
  $\breve{f}_{\varphi}: {\cal X}_{\varphi}\rightarrow {\cal Y}$
  of ${\cal O}_{\bullet}^{\ast}$-gerbes
   via the induced scheme-morphism $\hat{f}_{\varphi}$
   on their atlases;
  $\breve{f}_{\varphi}$
   is independent of the choices in the above discussion.
\end{lemma}

\smallskip

\noindent
It follows that
 the morphism $\varphi:(X^{A\!z},{\cal E})\rightarrow (Y,\alpha_B)$
 can be presented also as a morphism
 $\breve{\varphi}:({\cal X}^{A\!z},{\cal F})\rightarrow {\cal Y}$
 from an Azumaya ${\cal O}_X^{\ast}$-gerbe with a fundamental module,
 specified by a pair
 $$
  ({\cal O}_{\cal X} \subset \breve{\cal A}_{\varphi}  \subset
     {\cal O}_{\cal X}^{\Azscriptsize}
   :=\Endsheaf_{{\cal O}_{\cal X}}({\cal F})\;,\;
   \breve{f}_{\varphi}: {\cal X}_{\varphi} \rightarrow {\cal Y})\,.
 $$
The $\alpha_B$-twisted sheaf $\varphi_{\ast}{\cal E}$ on $Y$
 presented on the $\alpha_B$-admissible \'{e}tale cover $q:V\rightarrow Y$
 of $Y$ is given then by
 $\breve{q}^{\ast}
   \breve{f}_{\varphi,\ast}(_{{\cal O}_{{\cal X}_{\varphi}}}{\cal F})$.
Define
 $\breve{\varphi}_{\ast}{\cal F}
  := \breve{f}_{\varphi,\ast}(_{{\cal O}_{{\cal X}_{\varphi}}}{\cal F})$
 and let $Z=\Image\varphi\subset Y$
Then, $\breve{\varphi}_{\ast}{\cal F}$
 is supported on a substack of ${\cal Y}$
 that is a ${\cal O}_Z^{\ast}$-gerbe over $Z$.

\bigskip

In summary:

\bigskip

 
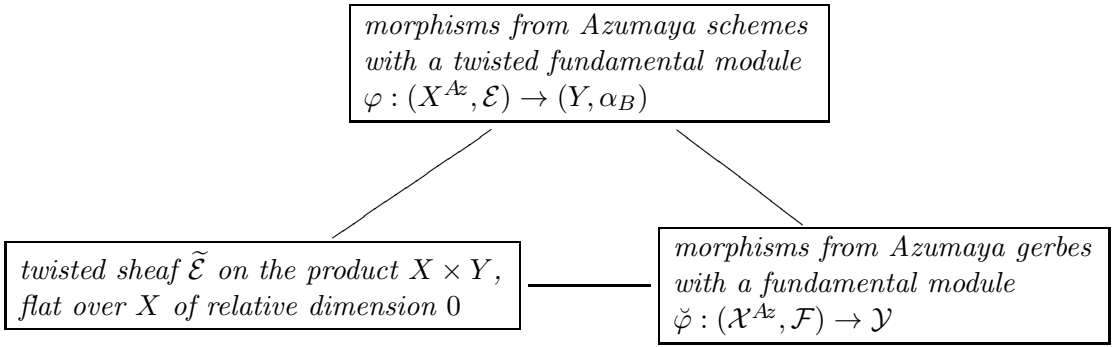
\xymatrix @R=8ex @C=-15ex {
  & \framebox[16.6em][c]{\parbox{15.6em}{\it
       morphisms from Azumaya schemes\\
       with a twisted fundamental module\\
       $\varphi:(X^{A\!z},{\cal E})\rightarrow (Y,\alpha_B)$}}
    \ar @{-} [ld]  \ar @{-}[rd]                  \\
  \framebox[17.8em][c]{\parbox{16.8em}{\it
   twisted sheaf $\widetilde{\cal E}$ on the product $X\times Y$,\\
   flat over $X$ of relative dimension $0$}}
   \ar @{-}[rr]
   && \framebox[15.6em][c]{\parbox{14.6em}{\it
       morphisms from Azumaya gerbes\\
       with a fundamental module\\
       $\breve{\varphi}:({\cal X}^{A\!z},{\cal F})\rightarrow {\cal Y}$}}
 } 

\bigskip
\bigskip

\noindent
Cf.~[L-L-S-Y: {\sc Figure}~2-2-1].

%
%
%

\bigskip

\section{The moduli stack of morphisms.}

\subsection{Family of D-branes in a {\boldmath $B$}-field background,
    twisted Hilbert polynomials, and boundedness.}

Some preparations toward the moduli problem of D-branes
 in a $B$-field background are given in this subsection.

\bigskip

\begin{flushleft}
{\bf Family of D-branes in a {\boldmath $B$}-field background.}
\end{flushleft}
The discussion in Sec.~2.2 applies also to a family.

\smallskip

\begin{definition}
{\bf [family of D-branes in {\boldmath $B$}-field background].}
{\rm
 Let
  $S$ be a base scheme$/{\Bbb C}$.
 An {\it $S$-family of morphisms}
  from Azumaya schemes with a fundamental module
  to $(Y,\alpha_B)$ consists of the following data:
  \begin{itemize}
   \item[$\cdot$]
    a flat family $X_S/S$ of schemes over $S$;

   \item[$\cdot$]
    a twisted coherent locally-free ${\cal O}_{X_S}$-module ${\cal E}_S$
     on $X_S/S$ of class
     $\alpha_S\in \check{C}^2_{\et}(X_S,{\cal O}_{X_S}^{\ast})$;

   \item[$\cdot$]
    a morphism
     $\,\varphi_S\,:\,
      (X_S^{A\!z}, {\cal E}_S)
       := (X_S,
           {\cal O}_{X_S}^{A\!z}:= \Endsheaf_{{\cal O}_{X_S}}({\cal E}_S),
           {\cal E}_S)\,
       \longrightarrow\, (Y,\alpha_B)\,$
     as defined in Definition~2.2.4.
  \end{itemize}
 Let
  $\,({\cal O}_{X_S} \subset {\cal A}_{\varphi_S}
               \subset {\cal O}_{X_S}^{A\!z}\;,\;
       f_{\varphi_S}: X_{\varphi_S}:=\boldSpec{\cal A}_{\varphi_S}
                  \rightarrow Y)\,$
   be the pair underlying $\varphi_S$  and
  $\pi_{\varphi_S}:X_{\varphi_S}\rightarrow X$ be the built-in morphism.

 Let
  $h:T\rightarrow S$ be a morphism of ${\Bbb C}$-schemes,
  $X_T = h^{\ast}X_S := T\times_S X_S$
    with the built-in $\hat{h}:X_T\rightarrow X_S$ that lifts $h$,
  ${\cal E}_T = \hat{h}^{\ast}{\cal E}_S$
   the pull-back $\hat{h}^{\ast}\alpha_S$($=:\alpha_T$)-twisted
    coherent locally-free ${\cal O}_{X_T}$-module,   and
  $(X_T^{A\!z},{\cal E}_T)
   = (X_T,
      {\cal O}_{X_T}^{A\!z}:=\Endsheaf_{{\cal O}_{X_T}}({\cal E}_T),
      {\cal E}_T)$.
 Then,
 the {\it pull-back} $h^{\ast}\varphi_S$ of $\varphi_S$ to $T$
  is the morphism
  $\varphi_T:(X_T^{A\!z},{\cal E}_T)\rightarrow (Y,\alpha_B)$
  with the underlying pair
   $$
    ({\cal O}_{X_T} \subset {\cal A}_{\varphi_T}
               \subset {\cal O}_{X_T}^{A\!z}\;,\;
       f_{\varphi_T}: X_{\varphi_T}:=\boldSpec{\cal A}_{\varphi_T}
                  \rightarrow Y)\,,
   $$
   where
   \begin{itemize}
    \item[$\cdot$]
     ${\cal A}_{\varphi_T}$ is
     the image ${\cal O}_{X_T}$-subalgebra of
       $\hat{h}^{\ast}{\cal A}_{\varphi_S}\rightarrow
        \hat{h}^{\ast}{\cal O}_{X_S}^{A\!z}={\cal O}_{X_T}^{A\!z}$,

    \item[$\cdot$]
     $f_{\varphi_T}$ is the composition of the morphisms
      $\xymatrix{\boldSpec{\cal A}_{\varphi_T} \ar@{^{(}->}[r]
        & h^{\ast}X_{\varphi_S}=\boldSpec(\hat{h}^{\ast}{\cal A}_S)
          \ar[r]^-{h^{\ast}f_{\varphi_S}} & Y}$.
   \end{itemize}
  Note that
   the minimal property of $X_{\varphi_T}$ is automatic  and
   the matching
    $\pi_{\varphi_T}^{\ast}\alpha_T=f_{\varphi_T}^{\ast}\alpha_B$
     follows from
      the matching
       $\pi_{\varphi_S}^{\ast}\alpha_S=f_{\varphi_S}^{\ast}\alpha_B$
        on $X_{\varphi_S}$ and
      the built-in inclusion
       $X_{\varphi_T}\hookrightarrow h^{\ast}X_{\varphi_S}$.
 In particular, let $\iota_s:s\rightarrow S$
  be a closed point of $S$.
 Then the {\it fiber} $\varphi_s$ of $\varphi_S$ over $s$
  is defined to be the morphism
  $\iota_s^{\ast}\varphi_S:
    (X_s^{A\!z},{\cal E}_s)\rightarrow (Y,\alpha_B)$.
}\end{definition}

\smallskip

\begin{remark}
{\it $[$surrogates in family$]$.} {\rm
 Note that in general $X_{\varphi_s}$, $s\in S$,
   do not form a flat family of schemes over $S$.
 See [L-L-S-Y: Remark 2.1.16] for more comments.
}\end{remark}

\smallskip

%
%

It follows from Lemma~2.2.7 that:

\smallskip

\begin{lemma}
{\bf [equivalent description
      via Azumaya-w/o-Azumaya-'n'-morphisms-w/o-morphisms].}
 An $S$-family of morphisms
  from Azumaya schemes with a fundamental module to $(Y,\alpha_B)$
  can be described equivalently by the following data:
  \begin{itemize}
   \item[$\cdot$]
    a flat family $X_S/S$ of schemes over $S$;

   \item[$\cdot$]
    a $\pr_1^{\ast}\alpha_S$-twisted
      coherent ${\cal O}_{X_S\times Y}$-module $\widetilde{\cal E}_S$
      on $X_S\times Y$
    that satisfies:
    \begin{itemize}
     \item[$(a)$]
      $\widetilde{\cal E}_S$ on $(X_S\times Y)/X_S$ is flat over $X_S$
       of relative dimension $0$ and of fixed relative length;

     \item[$(b)$]
      $\pr_1^{\ast}\alpha_S=\pr_2^{\ast}\alpha_B$
       on $\Supp\widetilde{\cal E}_S$.
    \end{itemize}
  \end{itemize}
 Here,
  $\pr_1:X_S\times Y\rightarrow X_S$  and
  $\pr_2:X_S\times Y\rightarrow Y$ are the projection maps.
\end{lemma}

\bigskip

\begin{flushleft}
{\bf Twisted Hilbert polynomials of a morphism.}
\end{flushleft}
\begin{definition-lemma}
{\bf [twisted Hilbert polynomial of sheaf].} {\rm (Cf.\ [Yo: Sec.~2.1].)
{\it
 Let
  $(W, {\cal O}_W(1))$ be a projective scheme
   with a fixed $\alpha$-twisted coherent locally-free
    ${\cal O}_W$-module ${\cal G}$
    in the class $\alpha\in \check{C}_{\et}(W,{\cal O}_W^{\ast})$
     with $[\alpha]\in \Br(W)$.
 Let ${\cal F}$ be an $\alpha$-twisted coherent ${\cal O}_W$-module.
 Then the function
  $$
   P_{{\cal G}, {\cal F}}:
    m \longmapsto
    \chi( ({\cal F}\otimes_{{\cal O}_W}{\cal G}^{\vee})(m) )\,,
   \hspace{2em} m\in {\Bbb Z}\,,
  $$
  is a polynomial in $m$ of degree $\dimm{\cal F}$.}
 We shall call it
  the {\it ${\cal G}$-twisted Hilbert polynomial} of ${\cal F}$
   on $(W,{\cal O}_W(1))$.
}\end{definition-lemma}

\smallskip

As tensoring by a twisted locally free sheaf leaves
 the flatness property of a twisted sheaf intact,
 one has the following proposition:

\smallskip

\begin{proposition}
 {\bf [invariance under flat deformation].} {\rm (Cf.\ [H-L] and [Ha].)}
 Let
  $S$ be a base scheme,
  $W_S\rightarrow S$ be a projective morphism
     with a relative ample line bundle ${\cal O}_{W_S/S}(1)$ on $W_S$,
  ${\cal G}_S$ be an $\alpha_S$-twisted coherent locally-free
   ${\cal O}_{W_S}$-module
   in the class $\alpha_S\in \check{C}_{\et}(W_S,{\cal O}_{W_S}^{\ast})$
       with $[\alpha_S]\in \Br(W_S)$.
 Let ${\cal F}_S$ be an $\alpha_S$-twisted coherent
  ${\cal O}_{W_S}$-module.
 Denote by ${\cal G}_s$ (resp.\ ${\cal F}_s$)
  the restriction of ${\cal G}_S$ (resp.\ ${\cal F}_S$)
  to the fiber $W_s$ of $W_S/S$ at a closed point $s\in S$.
 Then, if ${\cal F}_S$ is flat over $S$,
  the twisted Hilbert polynomial $P_{{\cal G}_s,{\cal F}_s}$
   is locally constant as a function of $s\in S$.
 When $S$ is reduced, the converse also holds.
\end{proposition}

\smallskip

Lemma~3.1.3
 motivates then the following definition:

\smallskip

\begin{definition}
{\bf [twisted Hilbert polynomial of morphism].} {\rm
 Assume that $Y$ is projective
  with a very ample line bundle ${\cal O}_Y(1)$.
 Fix an $\alpha$-twisted locally free coherent
  ${\cal O}_X$-module ${\cal G}$ on $(X, {\cal O}_X(1))$
  in the class $\alpha$ with $[\alpha]\in \Br(X)$.
 Let
  $\varphi:(X^{A\!z},{\cal E})\rightarrow (Y,\alpha_B)$
   be a morphism in the class $\alpha$  and
  $\widetilde{\cal E}$ be the $\pr_1^{\ast}\alpha$-twisted sheaf
   on $X\times Y$ that represents $\varphi$.
 (Here, $\pr_1:X\times Y\rightarrow X$ and
  $\pr_2:X\times Y\rightarrow Y$ are the projection maps.)
 Then, the {\it ${\cal G}$-twisted Hilbert polynomial
  $P_{{\cal G},\varphi}$ of $\varphi$}
  is defined to be $P_{pr_1^{\ast}{\cal G},\widetilde{\cal E}}$,
  where ${\cal O}_{X\times Y}(1)$
   is taken to be
   ${\cal O}_X(1)\boxtimes{\cal O}_Y(1)
    := \pr_1^{\ast}{\cal O}_X(1)\otimes\pr_2^{\ast}{\cal O}_Y(1)$.
}\end{definition}

\smallskip

\begin{lemma}
{\bf [invariance in a family].}
 Fix an $\alpha_S$-twisted locally free coherent
  ${\cal O}_{X_S}$-module ${\cal G}_S$ on $(X_S, {\cal O}_{X_S/S}(1))$
  in the class $\alpha_S$ with $[\alpha_S]\in \Br(X_S)$.
 Let $\varphi_S:(X_S^{A\!z},{\cal E}_S)\rightarrow (Y,\alpha_B)$
  be an $S$-family of morphisms from Azumaya schemes with
  a fundamental module to a projective $(Y,\alpha_B)$
  in the class $\alpha_S$.
 Then,
  the ${\cal G}_s$-twisted Hilbert polynomial $P_{{\cal G}_s,\varphi_s}$
   of $\varphi_s$ is locally constant as a function of $s\in S$.
\end{lemma}

\smallskip

\noindent
This is a consequence of Proposition~3.1.5.
Thus, a twisted Hilbert polynomial of a morphism gives the notation
 of {\it combinatorial type}
 of a morphism $\varphi:(X^{A\!z},{\cal E})\rightarrow (Y,\alpha_B)$.

\bigskip

\begin{flushleft}
{\bf Boundedness.}
\end{flushleft}
Recall first the following theorem from C\u{a}ld\u{a}raru [C\u{a}]:

\smallskip

\begin{theorem}
{\bf [equivalence of $\ModCategory(W,\alpha)$
      and $\ModCategory\mbox{-}{\cal A}$].}
 {\rm ([C\u{a}: Theorem 1.3.7].)}
 Let
  $\alpha\in \check{C}^2_{\et}(W,{\cal O}_W^{\ast})$
   with $[\alpha]\in \Br(W)$,
  ${\cal A}$ be an Azumaya algebra on $W$ with $[{\cal A}]=[\alpha]$,
  $\ModCategory(W,\alpha)$ be
   the category of ${\alpha}$-twisted ${\cal O}_W$-modules, and
  $\ModCategory\mbox{-}{\cal A}$
   be the category of right ${\cal A}$-modules on $W$.
 Let ${\cal G}$ be an $\alpha$-twisted coherent locally-free
  ${\cal O}_W$-module on $W$
  such that
   ${\cal A}\simeq \Endsheaf_{{\cal O}_W}({\cal G})$
    ($\simeq {\cal G}\otimes_{{\cal O}_W} {\cal G}^{\vee}$ canonically)
   as ${\cal O}_W$-algebras.
 Note that ${\cal G}$ is naturally a left ${\cal A}$-module.
 Then, the following pair of functors defines
  an equivalence of categories:
  $$
   \begin{array}{cccl}
    \ModCategory(W,\alpha)
     & \xymatrix{\ar@<.4ex>[rr] && \ar@<.4ex>[ll]}
     & \ModCategory\mbox{-}{\cal A} \\[.6ex]
    \bullet
     & \xymatrix{\ar@{|->}[rr] &&}
     & \bullet \otimes_{{\cal O}_W}{\cal G}^{\vee}\\[.2ex]
    \bullet\otimes_{\cal A}{\cal G}
     & \xymatrix{ && \ar@{|->}[ll]}
     & \bullet   &.
   \end{array}
  $$
\end{theorem}

\smallskip

The following proposition generalizes [L-L-S-Y: Proposition~2.3.1]:

\smallskip

\begin{proposition}
{\bf [boundedness of morphisms].}
 Assume that $Y$ is projective
  with a very ample line bundle ${\cal O}_Y(1)$.
 Let
  $(X_S/S, {\cal O}_{X_S/S}(1))$ be a flat family of projective schemes,
  ${\cal G}_S$ be an $\alpha$-twisted locally free coherent
   ${\cal O}_X$-module on $X_S$
    in the class $\alpha_S$ with $[\alpha_S]\in \Br(X_S)$,
  $(X^{A\!z}_S,{\cal E}_S)/S$ be a flat family of Azumaya schemes
   with a fundamental module over $X_S/S$ in the class $\alpha_S$,
  and $P$ be a polynomial in one variable.
 Then
  the set $\{\varphi_{\bullet}\}_{\bullet}$ of morphisms
    from fibers $(X^{A\!z}_s,{\cal E}_s)$ of $(X_S^{A\!z},{\cal E}_S)/S$
    to $(Y,\alpha_B)$
    with ${\cal G}_s$-twisted Hilbert polynomial
     $P_{{\cal G}_s, \varphi_s}=P$
   is bounded.
\end{proposition}

\begin{proof}
 Let $\pr_1$ be the projection map $X_s\times Y\rightarrow X_s$
   for $s\in S$.
 Observe that
 for the $\pr_1^{\ast}\alpha_s$-twisted ${\cal O}_{X_s\times Y}$-module
  $\widetilde{\cal E}_s$ that represents
  a morphism $\varphi_s:(X_s^{A\!z},{\cal E}_s)\rightarrow (Y,\alpha_B)$,
 there is an surjective homomorphism
  $\pr_1^{\ast}{\cal E}_s\rightarrow \widetilde{\cal E}_s$
  of $\pr_1^{\ast}\alpha_s$-twisted modules.
 As $\pr_1^{\ast}{\cal G}_s^{\vee}$ is locally free,
  one has the following exact sequence of the underlying
  ${\cal O}_{X_s\times Y}$-modules of
  the right
  $\Endsheaf_{{\cal O}_{X_s\times Y}}(\pr_1^{\ast}{\cal G}_s)$-modules
  in question:
  $$
   \pr_1^{\ast}{\cal E}_s
     \otimes_{{\cal O}_{X_s\times Y}}\pr_1^{\ast}{\cal G}_s^{\vee}\;
    \longrightarrow\;
    \widetilde{\cal E}_s
      \otimes_{{\cal O}_{X_s\times Y}}\pr_1^{\ast}{\cal G}_s^{\vee}\;
    \longrightarrow\; 0\,.
  $$
 The proposition follows now from Theorem~3.1.8 and
  [H-L: Lemma~1.7.6],
   which says that
    a family $\{F_i\}_{i\in I}$ of (ordinary) coherent sheaves
     on a projective scheme is bounded
    if and only if
      the set of Hilbert polynomials $\{P(F_i)\}_{i\in I}$ is finite
       and
      there is a coherent sheaf $F$ such that all $F_i$ admits
       surjective homomorphisms $F\rightarrow F_i$.

\end{proof}

%
%
%
%
%
%

\bigskip

\subsection{${\frak M}_{A\!z(X_S/S,\alpha_S)^f}(Y,\alpha_B)$
            is algebraic.}

Let
 $X_S/S$ be a (fixed) flat family of projective schemes over $S$,
 $\alpha_S\in \check{C}^2_{\et}(X_S,{\cal O}_{X_S}^{\ast})$
  with $[\alpha]\in \Br(X_S)$,  and
 ${\frak M}_{A\!z(X_S/S,\alpha_S)^f}(Y,\alpha_B)$
  be the moduli stack of morphisms
  from (non-fixed) Azumaya schemes with a fundamental module
   on (fixed) fibers $X_s$ of $X_S/S$
   to a (fixed) projective $(Y,\alpha_B)$ in the class $\alpha_s$.
As a sheaf of groupoids on the category $\Scheme_S$ of schemes over $S$
 with the fppf topology,
 $$
  {\frak M}_{A\!z(X_S/S,\alpha_S)^f}(Y,\alpha_B)(T)\;
  =\;\{ \varphi_T: (X_T^{A\!z}, {\cal E}_T)\rightarrow (Y,\alpha_B)\}\,,
 $$
 for $T\in \Scheme_S$.
Here,
 $X_T= T\times_S X_S$ with a built-in $X_T\rightarrow X_S$,
 $\alpha_T\in \check{C}^2_{\et}(X_T)$ the pull-back of $\alpha_S$,  and
 ${\cal E}_T$ is a $\alpha_T$-twisted
    coherent locally free ${\cal O}_{X_T}$-module.
A morphism $\varphi_{T,1}\rightarrow \varphi_{T,2}$
 for
 $$
  (\varphi_{T,1}: (X_{T,1}^{A\!z},{\cal E}_{T.1})\rightarrow (Y,\alpha_B))
  \;,
  (\varphi_{T,1}: (X_{T,1}^{A\!z},{\cal E}_{T.1})\rightarrow (Y,\alpha_B))\;
  \in\; {\frak M}_{A\!z(X_S/S,\alpha_S)^f}(Y,\alpha_B)(T)
 $$
 is an isomorphism
 $h:{\cal E}_{T,1}\stackrel{\sim}{\rightarrow}{\cal E}_{T,2}$
  of $\alpha_T$-twisted ${\cal O}_{X_T}$-modules
 that satisfies
  $\tilde{h}({\cal A}_{\varphi_{T,1}})={\cal A}_{\varphi_{T,2}}$  and
  the diagram
   \begin{eqnarray*}
    \xymatrix{
     X_{\varphi_{T,1}}\ar[rr]^{f_{\varphi_{T,1}}}     && Y\; .\\
     X_{\varphi_{T,2}}\ar[urr]_{f_{\varphi_{T,2}}}\ar[u]^{\widehat{h}}
                                                      &&
    }
   \end{eqnarray*}
   commutes.
 Here,
  $\tilde{h}: {\cal O}_{X_{T,1}}^{A\!z}
              \stackrel{\sim}{\rightarrow} {\cal O}_{X_{T,2}}^{A\!z}$
   and
  $\widehat{h}: X_{\varphi_{T,2}}
             \stackrel{\sim}{\rightarrow} X_{\varphi_{T,1}}$
  are $h$-induced isomorphisms
   of ${\cal O}_{X_T}$-algebras and $X_T$-schemes respectively.
The goal of this subsection is to prove
 that ${\frak M}_{A\!z(X_S/S,\alpha_S)^f}(Y,\alpha_B)$
  is an algebraic stack, locally of finite type.

Denote by
 ${\frak T}((X_S\times Y)/S, \pr_2^{\ast}\alpha_B)^{0/X_S/S}$
 the moduli stack of $\pr_2^{\ast}\alpha_B$-twisted coherent
  ${\cal O}_{X_s\times Y}$-modules on $X_s\times Y$, $s\in S$,
 that are flat over $X_s$ of relative dimension $0$.

Recall first the following proposition of Lieblich [Lie1]:
 (in the special case of schemes and
  with the sheaf {\boldmath $\mu$}$_r$
                          of groups of $r$-th roots of unity
  replaced by ${\cal O}_{\bullet}^{\ast}$)

\smallskip

\begin{proposition}
{\bf [stack of twisted coherent sheaves algebraic].}
{\rm ([Lie1: Proposition 4.1.1.1]).}
 Let
  $W$ be a scheme/$\,{\Bbb C}$,
  $\alpha\in \check{C}^2_{\et}(W)$,  and
  ${\frak T}(W,\alpha)$ be the stack of
   $\alpha$-twisted coherent ${\cal O}_W$-modules.
 Then,
  ${\frak T}(W,\alpha)$ is algebraic, locally of finite type.
\end{proposition}

\noindent
The proof follows from Artin's criteria for algebraic stacks,
 [Art] and [Schl].

\smallskip

\begin{proposition}
{\bf [${\frak T}((X_S\times Y)/S, \pr_2^{\ast}\alpha_B)^{0/X_S/S}$
      algebraic].}
 ${\frak T}((X_S\times Y)/S, \pr_2^{\ast}\alpha_B)^{0/X_S/S}$
 is an algebraic stack, locally of finite type.
\end{proposition}

\begin{proof}
 This follows from
  Proposition~3.2.1 and
  the observation that
   ${\frak T}((X_S\times Y)/S, \pr_2^{\ast}\alpha_B)^{0/X_S/S}$
   can be identified with the stack of morphisms
   from the fibers of the fixed $X_S/S$ to the connected components
   of ${\frak T}(Y,\alpha_B)$ that parameterizes $0$-dimensional
   $\alpha_B$-twisted ${\cal O}_Y$-modules.

\end{proof}

\smallskip

\begin{corollary}
{\bf [${\frak M}_{A\!z(X_S/S,\alpha_S)^f}(Y,\alpha_B)$ algebraic].}
 ${\frak M}_{A\!z(X_S/S,\alpha_S)^f}(Y,\alpha_B)$
  is an algebraic stack, locally of finite type.
\end{corollary}

\begin{proof}
 Let
  $\pr_1: X_S\times Y\rightarrow X_S$ and
   $\pr_2:X_S\times Y\rightarrow Y$ be projective maps  and
  $U\rightarrow X_S\times Y$ be an \'{e}tale cover of $X_S\times Y$
   that is both $\pr_1^{\ast}{\alpha_S}$-
    and $\pr_2^{\ast}\alpha_B$-admissible.
 Then the pair
  $(\pr_1^{\ast}\alpha_S,\pr_2^{\ast}\alpha_B)
    \in \check{C}^2_{\et}(X_S\times Y, {\cal O}_{X_S\times Y}^{\ast})
         \times
         \check{C}^2_{\et}(X_S\times Y, {\cal O}_{X_S\times Y}^{\ast})$
  of $2$-cocycles with value in ${\cal O}_{X_S\times Y}^{\ast}$
  determines the $2$-cochain
   $\pr_1^{\ast}\alpha_S-\pr_2^{\ast}\alpha_B$
   of ideal sheaves\footnote{In
                             this language, we have to allow
                             the term ``ideal sheaf" to include
                             also the nonproper one, i.e.\
                             ${\cal O}_{\bullet}$ itself.}
   of ${\cal O}_{X_S\times Y}$,
  and hence a $2$-cochain $(Z_{ijk})_{ijk}$ of closed subscheme,
   presented via $U$.
 The matching condition
   $\pi_{\varphi_s}^{\ast}\alpha_s=f_{\varphi_s}^{\ast}\alpha_B$
  of twists on $X_{s,\varphi_s}$
  for a $\varphi_s:(X_s^{A\!z},{\cal E}_s)\rightarrow (Y,\alpha_B)$
  is equivalent to
  the condition that
   \begin{itemize}
    \item[$(\dag)$]
     {\it
      the lift of the closed subscheme $\Supp\widetilde{\cal E}_s$
      of $X_S\times Y$ to $U\times_{X_S\times Y}U\times_{X_S\times Y}U$
      is contained in the $2$-cochain $(Z_{ijk})_{ijk}$}
      of closed subschemes
   \end{itemize}
  on the corresponding
   $\pr_1^{\ast}\alpha_s$-twisted sheaf $\widetilde{\cal E}_s$
   on $X_S\times Y$ that represents $\varphi_s$.
 Note that each $Z_{ijk}$ projects to a constructible subset
  of $X_S\times Y$ under the \'{e}tale morphism
  $U\times_{X_S\times Y}U\times_{X_S\times Y}U\rightarrow X_S\times Y$.
 Note also that,
  given a $\widetilde{\cal E}_s$ that represents a $\varphi_s$,
  it follows from
   the projectivity of $X_S/S$ and $Y$ and
   a generalization of the construction in [H-L: Sec.~2.2 and Chap.~4]
    to twisted sheaves,
    e.g.\ [H-S: Sec.~2], [Lie1: Sec.~4.1], [Yo: Sec.~2],
  that a small enough local chart $T$ of
   ${\frak T}((X_S\times Y)/S,\pr_2^{\ast}\alpha_B)^{0/X_S/S}$
   around $[\widetilde{\cal E}_s]$
   comes from a Quot-scheme construction
   on a pair of $\pr_2^{\ast}\alpha_B$-twisted locally free
   ${\cal O}_{X_S\times Y}$-modules on $X_S\times Y$
   for a local chart of $[\widetilde{\cal E}_s]$
   in ${\frak T}(X_S\times Y,\pr_2^{\ast}\alpha_B)$.
 In other words, $T$ is realized as the subscheme of
  the variety of homomorphisms
   ${\cal F}_1\stackrel{h}{\rightarrow} {\cal F}_0$
   for two fixed $\pr_2^{\ast}\alpha_B$-twisted locally free
   ${\cal O}_{X_S\times Y}$-modules,
   whose $\coker h$ correspond to
   objects of ${\frak T}(X_S\times Y,\pr_2^{\ast}\alpha_B)$.
 Condition~$(\dag)$ imposes now a system of determinantal-type
  constructible-subset conditions on $T$.
 Thus,
 it selects a local chart $T^{\prime}$ of
   ${\frak M}_{A\!z(X_S/S,\alpha_S)^f}(Y,\alpha_B)$ around
    $[\widetilde{\cal E}_s]$
  as the intersection of a system of determinantal-type
   constructible subset of $T$.
 Functoriality of the construction realizes
   ${\frak M}_{A\!z(X_S/S,\alpha_S)^f}(Y,\alpha_B)$
  then as a constructible substack of
  ${\frak T}((X_S\times Y)/S, \pr_2^{\ast}\alpha_B)^{0/X_S/S}$.
 The proposition now follows from Proposition~3.2.2.

\end{proof}

\smallskip

The following discussion shows that
 the twist-matching condition is a closed condition
 from the nature of the basic lemma below, which is immediate:

\smallskip

\begin{lemma} {\bf [basic].}
 Let
  $T$ be a discrete valuation ring with the field of fraction $K$,
  $R$ be a $T$-algebra,
  $M$ be an $R$-module that is flat over $T$, and
  $M_K:=M\otimes_T K$.
 Let $r\in R$ such that $r\cdot M_K=0$.
 Then, $r\cdot M=0$.
\end{lemma}

\smallskip

\begin{lemma}
{\bf [matching of twists as closed condition].}
 Let
  $T=\Spec R$, $R$ a discrete valuation ring,
   with the generic point $\eta$ and the closed point ${\mathbf 0}$,
  $h:T\rightarrow S$,
  $X_T=h^{\ast}X_S=T\times_SX_S$
   with the built-in map $\hat{h}:X_T\rightarrow X_S$,
  $\alpha_T=\hat{h}^{\ast}\alpha_S\,$;
  $\,\pr_1:X_T\times Y\rightarrow X_T$ and
    $\pr_2:X_T\times Y\rightarrow Y$ be the projection maps, and
  $\pi:X_T\times Y\rightarrow T$ be the built-in morphism.
 Let ${\cal F}_T$ be a $\pr_2^{\ast}\alpha_B$-twisted
  coherent ${\cal O}_{X_T\times Y}$-module on $X_T\times Y$
  such that
  \begin{itemize}
   \item[$\cdot$]
    ${\cal F}_T$ is flat over $T$,

   \item[$\cdot$]
    $\pr_1^{\ast}\alpha_{\eta}=\pr_2^{\ast}\alpha_B$
    on $\Supp({\cal F}_T|_{\eta})$.
  \end{itemize}
 Then,
  $\pr_1^{\ast}\alpha_T=\pr_2^{\ast}\alpha_B$ holds
  on $\Supp{\cal F}_T$ over $T$.
 The same statement holds also for ${\cal F}_T$
  being a $\pr_1^{\ast}\alpha_T$-twisted coherent
  ${\cal O}_{X_T\times Y}$-module on $X_T\times Y$.
\end{lemma}

\begin{proof}
 Let
  $p:U\rightarrow X_T\times Y$ be an \'{e}tale cover of $X_T\times Y$
   that is admissible to
   both $\pr_1^{\ast}\alpha_T$ and $\pr_2^{\ast}\alpha_B$  and
  $p^{(2)}: U^{(2)}:=U\times_{X_T\times Y}U\times_{X_T\times Y}U
                \rightarrow X_T\times Y$
   be the built-in morphism from the fibered product.
 Then,
  the pullback family $p^{(2),\ast}{\cal F}_T$ on $U^{(2)}$
  is flat over ${\bf 0}\in T$.
 The pair $(\pr_1^{\ast}\alpha_T,\pr_2^{\ast}\alpha_B)$ on $U^{(2)}$,
  each of which takes values in ${\cal O}_{U^{(2)}}^{\ast}$,
  determines the principal ideal sheaf\footnote{Here,
                              we allow a local generator of
                               a principal ideal sheaf
                               to be invertible,
                              cf.\ footnote~13.}
   $\pr_1^{\ast}\alpha_T - \pr_2^{\ast}\alpha_B$ of ${\cal O}_{U^{(2)}}$.
 The lemma now follows from Lemma~3.2.4.

\end{proof}

\smallskip

\begin{corollary}
{\bf [closed substack].}
 ${\frak M}_{A\!z(X_S/S,\alpha_S)^f}(Y,\alpha_B)$
  is a closed substack of
 ${\frak T}((X_S\times Y)/S,\pr_2^{\ast}\alpha_B)^{0/X_S/S}$.
\end{corollary}

\smallskip

\noindent
We remark that in general
 ${\frak T}((X_S\times Y)/S,\pr_2^{\ast}\alpha_B)^{0/X_S/S}$
 and, hence, ${\frak M}_{A\!z(X_S/S,\alpha_S)^f}(Y,\alpha_B)$
 are not closed in ${\frak T}(X_S\times Y,\pr_2^{\ast}\alpha_B)$.

\bigskip

\section{The case of holomorphic D-strings.}

In this section\footnote{The
                         current section continues
                         the previous work [L-L-S-Y] (D(2))
                         with Si Li and Ruifang Song, fall 2008.
                        C.-H.L.\ thank
                         them for the participation of
                          the biweekly Saturday D-brane working seminar,
                          spring 2008,  and
                         Liang Kong and S.L.\
                          for a communication/conversation
                          on further issues
                          while editing the current manuscript.},
 we consider the case
 when $X$ is a nodal/prestable curve and
      $Y$ is a smooth projective variety.

\bigskip

\subsection{The moduli stack
    ${\frak M}_{A\!z(g,r,\chi)^f}(Y,\alpha_B;\beta)$
    of morphisms from Azumaya prestable curves to $(Y,\alpha_B)$.}

\begin{flushleft}
{\bf The big vs.\ the small moduli problem.}
\end{flushleft}
Recall the following lemma,
 which follows from the normalization sequence and Tsen's Theorem:
 (see [Lie1: Sec.~5.1.1] for more general discussions.)

\smallskip

\begin{lemma}
{\bf [$\Br(C)$ vanishes].} {\rm (Cf.\ [Lie1: Lemma~5.1.1.1].)}
 Let $C$ be a nodal curve. Then $\Br(C)=0$.
\end{lemma}

\smallskip

\noindent
Thus,
for
 $C$ a prestable curve and
 $\alpha\in\check{C}_{\et}(C,{\cal O}_C^{\ast})$
   with $[\alpha]\in \Br(C)$,
$[\alpha]$ indeed vanishes and
 \begin{itemize}
  \item[$\cdot$]
   there is an $\alpha$-twisted line bundle ${\cal L}$ on $C$
    and
   the correspondence
    $$
     \begin{array}{ccc}
      \ModCategory(C,\alpha)
       & \xymatrix{\ar[rr] &&}        & \ModCategory(C) \\[.6ex]
      {\cal F}
       & \xymatrix{\ar@{|->}[rr] &&}
       & {\cal F}\otimes_{{\cal O}_C}{\cal L}^{\vee}
     \end{array}
    $$
   is an equivalence of categories,

  \item[$\cdot$]
   any Azumaya algebra ${\cal A}$ over $C$ is isomorphic to
    $\Endsheaf_{{\cal O}_C}(\cal E)$
    for some (ordinary) locally free ${\cal O}_C$-module ${\cal E}$.
 \end{itemize}

Despite the fact that these special features reduce the study
 of Azumaya schemes and modules on prestable curve to the case
 as in [L-L-S-Y],
the moduli stacks ${\frak M}_{\alpha}(Y,\alpha_B)$
 of morphisms from Azumaya prestable curves with
 an $\alpha$-twisted fundamental module to $(Y,\alpha_B)$
 in general are {\it not} isomorphic for different choices
 of $\alpha$'s.
The image of such morphisms for different $\alpha$
  are in general distinct
 due to the effect of $\alpha_B$ via the twist-matching condition.
Thus, one has two moduli problems:
 \begin{itemize}
  \item[(1)]
  \parbox[t]{13.6em}{{\it The big moduli problem}$\,$:}
   moduli of morphisms from Azumaya prestable curves with
    a possibly twisted fundamental module to $(Y,\alpha_B)$.

  \item[(2)]
  \parbox[t]{13.6em}{{\it The small moduli problem}$\,$:}
   moduli of morphisms from Azumaya prestable curves with
     an ordinary/untwisted fundamental module to $(Y,\alpha_B)$.
 \end{itemize}
While Problem (1) is a final goal,
 in this work we address only Problem (2),
 a sub-problem to Problem (1).

\bigskip

\begin{flushleft}
{\bf The small moduli problem.\footnote{This
                       theme is taken/adapted from [L-L-S-Y].
                      Readers are referred ibidem for more details.}}
\end{flushleft}
Let
 ${\frak M}_{A\!z(g,r,\chi)^f}(Y,\alpha_B;\beta)$
 be the moduli stack of morphisms
 $\varphi:(C^{A\!z},{\cal E})\rightarrow (Y,\alpha_B)$
 from (unfixed) Azumaya prestable curves with a fundamental module
 to $(Y,\alpha)$ (cf.\ Definition~2.2.4)
 of type $(g,r,\chi,|\,\beta)$ in the sense that:
 \begin{itemize}
  \item[$\cdot$]
   $C$ has (arithmetic) genus $g$,

  \item[$\cdot$]
   ${\cal E}$ has rank $r$ and
   Euler characteristic $\chi=\degree{\cal E} + r(1-g)$,

  \item[$\cdot$]
   the image curve class $\varphi_{\ast}[C]=\beta\in N_1(Y)$.
 \end{itemize}
([L-L-S-Y: Definition~2.1.11, Definition~2.1.12,
           Definition~2.1.13, Remark~2.1.14, and Lemma~2.2.4].)
Then

\smallskip

\begin{proposition}
{\bf [${\frak M}_{A\!z(g,r,\chi)^f}(Y,\alpha_B;\beta)$ algebraic].}
 ${\frak M}_{A\!z(g,r,\chi)^f}(Y,\alpha_B;\beta)$
 is an algebraic stack, locally of finite type.
\end{proposition}

\smallskip

\noindent
This follows from [L-L-S-Y: Sec.~3.2] and Corollary~3.2.3.

\bigskip

\subsection{Fillability/valuation-criterion property of
            ${\frak M}_{A\!z(g,r,\chi)^f}(Y,\alpha_B;\beta)$.}

We prove in this subsection
 the following fillability/valuation-criterion property of
 the moduli stack ${\frak M}_{A\!z(g,r,\chi)^f}(Y,\alpha_B;\beta)$.
This indicates that ${\frak M}_{A\!z(g,r,\chi)^f}(Y,\alpha_B;\beta)$
 is a sufficiently large moduli space
 for the study of more restrictive moduli problems for curves.

\smallskip

\begin{proposition}
{\bf [valuation-criterion property].}
 Let
  $S=\Spec R$, $R$ a discrete valuation ring,
   with the generic point $\eta$ and the closed point ${\mathbf 0}$
   and
  $f:\eta\rightarrow {\frak M}_{A\!z(g,r,\chi)^f}(Y,\alpha_B;\beta)$
   be a morphism.
 Then, after a base change on $S$ if necessary,
  $f$ extends to a morphism
  $\hat{f}:S\rightarrow {\frak M}_{A\!z(g,r,\chi)^f}(Y,\alpha_B;\beta)$.
\end{proposition}

\smallskip

\noindent
(However, the extension in general is not unique.)

\smallskip

\begin{proof}
 Let
  $C_{\eta}$ be a flat family of prestable curves
   (of genus $g$) over $\eta$ and
  $\widetilde{\cal E}_{\eta}$
   be the coherent ${\cal O}_{C_{\eta}\times Y}$-module
   on $C_{\eta}\times Y$,
    flat over $C_{\eta}$
     with relative dimension $0$ (and of relative length $r$),
  that corresponds to the morphism
  $f:\eta \rightarrow {\frak M}_{A\!z(g,r,\chi)^f}(Y,\alpha_B;\beta)$.
 Let
  $\pr_1:C_{\eta}\times Y\rightarrow C_S$
   and $\pr_2:C_{\eta}\times Y\rightarrow Y$ be the projection maps.
 {From} the projection formula,
  $\pr_{1,\ast}({\cal F}
     \otimes_{{\cal O}_{C_{\eta}\times Y}}\pr_1^{\ast}{\cal L})
   \simeq (\pr_{1,\ast}{\cal F})\otimes_{{\cal O}_{C_{\eta}}}{\cal L}$
  for an ${\cal O}_{C_{\eta}\times Y}$-module ${\cal F}$ and
      a coherent locally free ${\cal O}_{C_{\eta}}$-module ${\cal L}$,
 we may assume,
   after tensoring $\pr_1^{\ast}{\cal L}$ for an appropriate
   relative ample line bundle ${\cal L}$ on $C_{\eta}/{\eta}$,
  that ${\cal E}_{\eta} := \pr_{1,\ast}\widetilde{\cal E}_{\eta}$
  fits into the exact sequence
   $$
    {\cal O}_{C_{\eta}}^{\;\oplus N}\;
     \longrightarrow\; {\cal E}_{\eta}\; \longrightarrow\; 0\,.
   $$
   for some $N\gg 0$.
 It follows from the built-in exact sequence
  $\pr_1^{\ast}{\cal E}_{\eta}
              \rightarrow \widetilde{\cal E}_{\eta}\rightarrow 0$
  that
  $$
   {\cal O}_{C_{\eta}\times Y}^{\;\oplus N}\;
    \longrightarrow\; \widetilde{\cal E}_{\eta}\;
    \longrightarrow\; 0\,.
  $$
 Regarding this as an exact sequence of
  ${\cal O}_{C_{\eta}\times Y}$-modules on
  $(C_{\eta}\times Y)/C_{\eta}/S$,
 one obtains a morphism
  $$
   \tilde{f}_{\eta}\; :\;
    C_{\eta}/\eta\;
    \longrightarrow\; \Quot_Y({\cal O}_Y^{\oplus N}, r)\,,
  $$
  where $\Quot_Y({\cal O}_Y^{\oplus N}, r)$ is the Quot-scheme
   that parameterizes the $0$-dimensional quotient sheaves
    of ${\cal O}_Y^{\oplus N}$ with length $r$.
 Rigidifying $\tilde{f}_{\eta}$
  as a prestable map from a curve over $\eta$ if necessary,
 the properness property of the moduli stack
   $\overline{\cal M}_{\bullet}
                     (\Quot_Y({\cal O}_Y^{\oplus N}, r),\,\bullet\,)$
   of stable maps to the projective scheme
   $\Quot_Y({\cal O}_Y^{\oplus N}, r)$
  implies that, subject to a base change\footnote{For
                        the simplicity of notation and expressions,
                        a base change on $S$ and the new family over $S$
                        will still be denoted by $S$ and $C_S/S$
                        respectively.}
  on $S$,
 $\tilde{f}_{\eta}$ extends to a morphism
  $\tilde{f}_S:C_S/S\rightarrow \Quot_Y({\cal O}_Y^{\oplus N}, r)$,
   where $C_S/S$ is a flat family of prestable curves that extends
   $C_{\eta}/\eta$ (after the above base change).
 The associated quotient sheaf $\widetilde{\cal E}_S$
  on $C_S\times Y$
  extends $\widetilde{\cal E}_{\eta}$ on $C_{\eta}\times Y$
  and has the property that $\widetilde{\cal E}_S$
  is flat over $C_S$ of relative dimension $0$ and relative length $r$.
 This defines thus
  $\hat{f}:S\rightarrow {\frak M}_{A\!z(g,r,\chi)^f}(Y,\beta)$.
 Corollary~3.2.6 implies then
  that $\hat{f}$ indeed has the image in
  ${\frak M}_{A\!z(g,r,\chi)^f}(Y,\alpha_B;\beta)$ and,
  hence, is an extension of
  $f:S-\{\mathbf 0\} \rightarrow
     {\frak M}_{A\!z(g,r,\chi)}^f(Y,\alpha_B;\beta)$,
  after a base change.
 This proves the proposition.

\end{proof}

\bigskip

\section{The extension by the sheaf ${\cal D}$ of differential operators.}

In this section, the second effect
  - namely, the deformation quantization
    of both the target space(-time) and D-brane world-volumes -
 of the background $B$-field to a smooth D-brane world-volume $X$
 along the line of the Polchinski-Grothendieck Ansatz
 is brought into consideration as well.
We focus on the case when the deformation quantizations that occur
 are modelled directly on that for the phase space in quantum mechanics
 and when the study in Sec.~2.2 can be extended/applied immediately.
The special case of morphisms from $X$ with the new structure
  to a target-space $Y$ being the total space {\boldmath $\Omega$}$_W$
  of the cotangent bundle $\Omega_W$ of a smooth variety $W$
 is considered.
An application of this gives the notion of
 deformation quantizations of the spectral curves that appear
 in Hitchin's integrable systems.
For language simplicity, we use the analytic topology
 for smooth varieties in the discussion below
 whenever it is more convenient.

\bigskip

\subsection{Azumaya schemes with a fundamental module
            with a flat connection.}

The discussion in Sec.~2.2 has a direct generalization to incorporate
 the sheaf ${\cal D}$ of algebras of differential operators
 and ${\cal D}$-modules.

\bigskip

\begin{flushleft}
{\bf Weyl algebras, the sheaf ${\cal D}$
     of differential operators, and ${\cal D}$-modules.}
\end{flushleft}
Let
 $X$ be a smooth variety over ${\Bbb C}$,
 $\Theta_X=\Dersheaf_{{\Bbb C}}({\cal O}_X,{\cal O}_X)$
  be the sheaf of ${\Bbb C}$-derivations on ${\cal O}_X$,  and
 $\Omega_X$ be the sheaf of K\"{a}hler differentials on $X$.
We recall a few necessary objects and facts for our study.
Their details are referred to
 [Be3], [Bj], and [B-E-G-H-K-M]$\,$:$\,$\footnote{See
  Bernstein [Be3: {\S}0.\ Introduction], Bj\"{o}rk [Bj: Introduction],
   and Borel [B-E-G-H-K-M: Introduction]
  for a list of people who contribute to the early development
  of the subject.}
\begin{itemize}
 \item[(1)]
  the {\it Weyl algebra}
  $$
    A_n({\Bbb C})\; :=\;
     {\Bbb C}\langle x_1,\,\cdots\,,\,x_n,
       \partial_1,\,\cdots\,,\,\partial_n\rangle
         /([x_i,x_j]\,,\, [\partial_i,\partial_j]\,,\,
           [\partial_i,x_j]-\delta_{ij}\, :\, 1\le i,j\le n)\,,
  $$
  which is the algebra of differential operators acting on
   ${\Bbb C}[x_1,\,\cdots\,,\,x_n]$ by formal differentiation;
  here,
    ${\Bbb C}\langle\,\cdots\,\rangle$ is
     the unital associative ${\Bbb C}$-algebra generated
      by elements $\cdots$ indicated,
    $[\;\,,\,\;]$ is the commutator,
    $\delta_{ij}$ is the Kronecker delta, and
    $(\,\cdots\,)$ is the $2$-sided ideal generated by $\cdots$ indicated;

 \item[(2)]
 the {\it sheaf ${\cal D}_X$
  of (linear algebraic) differential operators}
  on $X$,
  which is the sheaf of unital associative algebras
  that extends ${\cal O}_X$ by new generators from the sheaf $\Theta_X$;

 \item[(3)]
  {\it ${\cal D}_X$-modules}
   (or directly {\it ${\cal D}$-modules} when $X$ is understood),
   which are sheaves on $X$ on which ${\cal D}_X$ acts from the left.
 %
 %
\end{itemize}

\smallskip

\begin{lemma}
{\bf [$A_n({\Bbb C})$ simple].}
 $A_n({\Bbb C})$ is a simple algebra:
  the only $2$-sided ideal therein is the zero ideal $(0)$.
\end{lemma}

\smallskip

\begin{proposition}
{\bf [${\cal O}$-coherent ${\cal D}$-module].}
 Let ${\cal M}$ be a ${\cal D}_X$-module
  that is coherent as an ${\cal O}_X$-module.
 Then, ${\cal M}$ is ${\cal O}_X$-locally-free.
 Furthermore, in this case,
 the action of ${\cal D}_X$ on ${\cal M}$ defines a flat connection
  $\nabla:{\cal M}\rightarrow {\cal M}\otimes\Omega_X$ on ${\cal M}$
  by assigning $\nabla_{\!\xi}\,s=\xi\cdot s$
   for $s\in{\cal M}$ and $\xi\in \Theta_X$;
 the converse also holds.
 This gives an equivalence of categories:
 $$
  \left\{\rule{0em}{1.2em}
   \begin{array}{c}
    \mbox{${\cal O}_X$-coherent ${\cal D}_X$-modules}
   \end{array}
  \right\}\;
  \longleftrightarrow\;
  \left\{
   \begin{array}{l}
     \mbox{coherent locally free ${\cal O}_X$-modules}\\
     \mbox{with a flat connection}
   \end{array}
  \right\}\,.
 $$
\end{proposition}

\bigskip

\begin{flushleft}
{\bf ${\cal D}$ as the structure sheaf of
     the deformation quantization of the cotangent bundle.}
\end{flushleft}
{From} the presentation of the Weyl algebra $A_n({\Bbb C})$,
 which resembles the quantization of a classical phase space
  with the position variable $(x_1,\,\cdots\,,\,x_n)$ and
  the dual momentum variable
   $(p_1,\,\cdots\,,\,p_n)=(\partial_1,\,\cdots\,,\,\partial_n)$,
  and
the fact that
 ${\cal D}_X$ is locally modelled on the pull-back of
 $A_n({\Bbb C})$ over ${\Bbb A}^n$ under an \'{e}tale morphism
 to ${\Bbb A}^n$,
the sheaf ${\cal D}_X$ of algebras with the built-in inclusion
 ${\cal O}_X\subset {\cal D}_X$ can be thought of
 as the structure sheaf of a noncommutative space
 from the quantization\footnote{The word
                         ``quantization" has received various meanings
                          in mathematics. Here, we mean solely the
                          one associated to quantum mechanics.
                          This particular quantization is also called
                           {\it deformation quantization}.}
 of the cotangent bundle, i.e.\ the total space
 $\mbox{\boldmath $\Omega$}_X$ of the sheaf $\Omega_X$, of $X$.

\smallskip

\begin{definition}
{\bf [canonical deformation quantization of cotangent bundle].}
{\rm
 We will formally denote this noncommutative space by
  $\Space {\cal D}_X =: Q\mbox{\boldmath $\Omega$}_X$ and
 call it
  the {\it canonical deformation quantization}
  of $\mbox{\boldmath $\Omega$}_X$.
}\end{definition}

\smallskip

A special class of {\it morphisms} from or to $\Space {\cal D}_X$
 can be defined contravariantly
 as homomorphisms of sheaves of ${\Bbb C}$-algebras.

\smallskip

\begin{example}
{\bf [$A_n({\Bbb C})$].} {\rm
 The noncommutative space $\Space (A_n({\Bbb C}))$
  defines a deformation quantization of
  $\mbox{\boldmath $\Omega$}_{{\Bbb A}^n}$.
 Recall the presentation of $A_n({\Bbb C})$.
 The ${\Bbb C}$-algebra homomorphism
  $$
   \begin{array}{cccccl}
    f_{(k)}^{\sharp} & :
     & {\Bbb C}[y_1,\,\cdots\,,\,y_n]
     & \longrightarrow     & A_n({\Bbb C}) \\[.6ex]
    && y_i & \longmapsto   & x_i\,,        & i=1,\,\ldots\,,\,k, \\[.6ex]
    && y_j & \longmapsto   & \partial_j\,, & j=k+1,\,\ldots\,,\, n\,,
   \end{array}
   $$
 defines a dominant morphism
  $f_{(k)}:\Space (A_n({\Bbb C}))\rightarrow {\Bbb A}^n$,
  $k=0,\,\ldots\,,\,n$.
 The ${\Bbb C}$-algebra automorphism
  $A_n({\Bbb C})\rightarrow A_n({\Bbb C})$ with
   $x_i\mapsto \partial_i$ and $\partial_i\mapsto -x_i$
   defines the {\it Fourier transform} on $\Space (A_n({\Bbb C}))$.
 Note that, since $A_n({\Bbb C})$ is simple,
  any morphisms to $\Space (A_n({\Bbb C}))$ is dominant
  (i.e.\ the related ${\Bbb C}$-algebra homomorphism
   from $A_n({\Bbb C})$ is injective).
}\end{example}

\bigskip

\begin{flushleft}
{\bf $\alpha$-twisted ${\cal O}_X$-coherent ${\cal D}_X$-modules and
     enlargements of ${\cal O}_X^{A\!z}$ by ${\cal D}_X$.}
\end{flushleft}
Let
 $\alpha\in \check{C}_{\et}(X,{\cal O}_X^{\ast})$  and
 $\,{\cal F}\; =\; ( \{U_i\}_{i\in I},\,
                   \{{\cal F}_i\}_{i\in I},\,
                   \{\phi_{ij}\}_{i,j\in I} )\,$
  be an $\alpha$-twisted ${\cal O}_X$-module.

\smallskip

\begin{definition}
{\bf [connection on ${\cal F}$].} {\rm
 A {\it connection} $\nabla$ on ${\cal F}$
  is a set $\{\nabla_i\}_{i\in I}$
  where
   $\nabla_i: {\cal F}_i \rightarrow
      {\cal F}_i\otimes_{{\cal O}_{U_i}}\Omega_{U_i}$
   is a connection on ${\cal F}_i$,
 that satisfies
  $\phi_{ij}\circ (\nabla_i|_{U_{ij}})
    = (\nabla_j|_{U_{ij}})\circ \phi_{ij}$.
 $\nabla$ is said to be {\it flat}
  if $\nabla_i$ is flat for all $i\in I$.
}\end{definition}

\smallskip

\noindent
Note that the existence of an $\alpha$-twisted ${\cal O}_X$-module
 with a connection imposes
 a condition on $\alpha$ that $\alpha$ has a presentation
 $(\alpha_{ijk})_{ijk}$ with
 $d\alpha := (d\alpha_{ijk})_{ijk}=(0)_{ijk}$;
    i.e.\ $\alpha_{ijk}\in{\Bbb C}^{\ast}$ for all $i,j,k$.

As the proof of Proposition~5.1.2 is local,
 it generalizes to
 $\alpha$-twisted ${\cal O}_X$-coherent ${\cal D}_X$-modules$\,$:

\smallskip

\begin{proposition}
{\bf [$\alpha$-twisted ${\cal O}$-coherent ${\cal D}$-module].}
 Let ${\cal M}$ be a ${\cal D}_X$-module
  that is $\alpha$-twisted ${\cal O}_X$-coherent.
 Then, ${\cal M}$ is an $\alpha$-twisted ${\cal O}_X$-locally-free.
 Furthermore, in this case,
 the action of ${\cal D}_X$ on ${\cal M}$ defines a flat connection
  $\nabla:{\cal M}\rightarrow {\cal M}\otimes\Omega_X$ on ${\cal M}$
  by assigning $\nabla_{\!\xi}\,s=\xi\cdot s$
   for $s\in{\cal M}$ and $\xi\in \Theta_X$;
 the converse also holds.
 This gives an equivalence of categories:
 $$
  \left\{\rule{0em}{1.2em}
   \begin{array}{c}
    \mbox{$\alpha$-twisted ${\cal O}_X$-coherent ${\cal D}_X$-modules}
   \end{array}
  \right\}\;
  \longleftrightarrow\;
  \left\{
   \begin{array}{l}
     \mbox{$\alpha$-twisted coherent locally free ${\cal O}_X$-}\\
     \mbox{modules with a flat connection}
   \end{array}
  \right\}\,.
 $$
\end{proposition}

\smallskip

Let ${\cal E}$ be an $\alpha$-twisted ${\cal O}_X$-coherent
  ${\cal D}_X$-module.
Then the ${\cal D}_X$-module structure on ${\cal E}$ induces
  a natural ${\cal D}_X$-module structure on the (ordinary)
  ${\cal O}_X$-module
  ${\cal O}_X^{A\!z} := \Endsheaf_{{\cal O}_X}({\cal E})$.
We will denote both the connection on ${\cal E}$ and
 on ${\cal O}_X^{A\!z}$ by $\nabla$.
As
 both ${\cal O}_X^{A\!z}:=\Endsheaf_{{\cal O}_X}({\cal E})$ and
  ${\cal D}_X$ act now on ${\cal E}$  and
 ${\cal D}_X$ acts also on ${\cal O}_X^{A\!z}$,
one can define a sheaf ${\cal O}_X^{A\!z,{\cal D}}$ of unital
 associative algebras generated by ${\cal O}_X^{A\!z}$ and ${\cal D}_X$
 as follows:
 \begin{itemize}
  \item[$\cdot$]
   Over a (Zariski) open subset $U$ of $X$,
    ${\cal O}_X^{A\!z,{\cal D}}(U)$ is the unital associative
    ${\Bbb C}$-algebra generated by
    ${\cal O}^{A\!z}_X(U)\cup {\cal D}_X(U)$
   subject to the following rules$\,$:
    \begin{itemize}
     \item[(1)]
      for $\phi_1,\,\phi_2\in {\cal O}_X^{A\!z}(U)$,
      $\phi_1\cdot\phi_2\in {\cal O}_X^{A\!z,{\cal D}}(U)$
       coincides with the existing
       $\phi_1\phi_2\in {\cal O}_X^{A\!z}(U)\,$;

     \item[(2)]
      for $\eta_1,\,\eta_2\in {\cal D}_X(U)$,
      $\eta_1\cdot\eta_2\in {\cal O}_X^{A\!z,{\cal D}}(U)$
       coincides with the existing
       $\eta_1\eta_2\in{\cal D}_X(U)\,$;

     \item[(3)] ({\it Leibniz rule})\hspace{1em}
      for $\phi\in {\cal O}_X^{A\!z}(U)$ and
           $\xi\in\Theta_X(U)\subset {\cal D}_X(U)$,
       $$
        \xi\cdot\phi\; =\; (\nabla_{\!\xi}\,\phi)\,+\, \phi\cdot \xi\,.
       $$
    \end{itemize}
 \end{itemize}
In notation,
 ${\cal O}_X^{A\!z,{\cal D}}
    := {\Bbb C}\langle {\cal O}_X^{A\!z},{\cal D}_X\rangle^{\nabla}$.

\smallskip

\begin{definition}
{\bf [Azumaya quantum scheme with fundamental module].} {\rm
The noncommutative space
 $$
  (X^{A\!z,{\cal D}},\,{\cal E}^{\nabla})\;
  :=\; (X,\,
       {\cal O}_X^{A\!z,{\cal D}}
       = {\Bbb C}\langle
          \Endsheaf_{{\cal O}_X}({\cal E}), {\cal D}_X\rangle^{\nabla},\,
       ({\cal E},\nabla) )
 $$
 will be called an {\it Azumaya quantum scheme with a fundamental module
 in the class $\alpha$}.
}\end{definition}

\smallskip

\noindent
Caution that ${\cal O}_X\subset {\cal O}_X^{A\!z,{\cal D}}$
 in general does not lie in the center of ${\cal O}_X^{A\!z,{\cal D}}$.

\smallskip

\begin{remark}
{$[\,$${\cal E}^{\nabla}$ as a module
      over $\Space({\cal O}_X^{A\!z,{\cal D}}) $$\,]$.}
{\rm
 The full notation for $X^{A\!z,{\cal D}}$ in Definition~5.1.7
   is meant to make two things manifest:
   \begin{itemize}
    \item[(1)]
     There is a built-in diagram of dominant morphisms of $X$-spaces$\,$:
      $$
       \xymatrix @R=1em @C=-1em {
        & X^{A\!z,{\cal D}} := \Space {\cal O}_X^{A\!z,{\cal D}}
          \ar[ld] \ar[rd] \ar[dd] & \\
        **[l]X^{A\!z} := \Space {\cal O}_X^{A\!z} \ar[rd]
         && **[r]Q\mbox{\boldmath $\Omega$}_X := \Space {\cal D}_X\,.
                 \ar[ld] \\
        & X  &
       }
      $$
     $\Space {\cal O}_X^{A\!z,{\cal D}}$ is the major space
        one should focus on.
     The other three spaces
       - $\Space {\cal O}_X^{A\!z}$, $\Space {\cal D}_X$, and $X$ -
      should be treated as auxiliary spaces that are built into
      the construction to encode a special treatment
      that takes care of the issue of
      localizations of noncommutative rings in the current situation;
      cf.\ the next item.

     \item[(2)]
      Despite the fact that ${\cal O}_X$ is in general not
             in the center of ${\cal O}_X^{A\!z,{\cal D}}$,
       there is a notion of localization and open sets on
       $\Space {\cal O}^{A\!z,{\cal D}}$ induced by those on $X$.
      I.e.\ $\Space{\cal O}_X^{A\!z,{\cal D}}$ has a built-in topology
       induced from the (Zariski) topology of $X$.
      Thus, one can still have the notion of
       {\it gluing systems of morphisms} and {\it sheaves}
       with respect to this topology.
   \end{itemize}
  In particular, ${\cal E}^{\nabla}$ is a sheaf of
   ${\cal O}_X^{A\!z,{\cal D}}$-modules
   supported on the whole $\Space{\cal O}^{A\!z,{\cal D}}$
   with this topology.
}\end{remark}

\smallskip

\begin{remark}
{$[\,$Azumaya algebra over ${\cal D}_X$$\,]$.}
{\rm
 Note that ${\cal O}_X^{A\!z,{\cal D}}$ can also be thought of
  as an {\it Azumaya algebra over ${\cal D}_X$}
  in the sense that it is a sheaf of algebras on $X$,
  locally modelled on the matrix ring
   $M_r({\cal D}_U)$ over ${\cal D}_U$
  for $U$ an affine \'{e}tale-open subset of $X$.
}\end{remark}

\smallskip

\begin{remark}
{\it $[\,$partially deformation-quantized target$\,]$.} {\rm
 {From} the fact that Weyl algebras are simple,
  it is anticipated that a morphism to
  a totally deformation-quantized space $Y=\mbox{\boldmath $\Omega$}_W$
  is a dominant morphism.
 In general, one may take $Y$ to be a partial deformation quantization
  of a space along a foliation.
 E.g.\  a deformation quantization of {\boldmath $\Omega$}$_{W/B}$
  along the fibers of a fibration $W/B$.
 For compact $Y$, one may consider
  the deformation quantization along torus fibers
  of a space fibered by even-dimensional tori.\footnote{Though
                                   we do not touch this here,
                                   readers should be aware that
                                    this is discussed in numerous
                                    literatures.}
 (Cf.~Example~5.1.11.)
}\end{remark}

\bigskip

\begin{flushleft}
{\bf Higgsing and un-Higgsing of D-branes via deformations of
morphisms.}
\end{flushleft}
Same as the situation studied in [L-Y1], [L-L-S-Y], [L-Y2], and [L-Y3],
the Higgsing and un-Higgsing of D-branes can occur
 when we deform morphisms in the current situation.

\smallskip

\begin{example}
{\bf [Higgsing/un-Higgsing of D-brane].} {\rm
 Let
  $(X^{A\!z,{\cal D}},{\cal E}^{\nabla})$
   be the affine Azumaya quantum scheme with a fundamental module
   associated to
    the ring
     $R:= {\Bbb C}\langle M_2({\Bbb C}[z]), \partial_z\rangle$
     (with the implicit relation $[\partial_z,z]=1$ and
           the identification of ${\Bbb C}[z]$
                   with the center of $M_2({\Bbb C}[z])$)  with
    the $R$-module $N:={\Bbb C}[z]\oplus {\Bbb C}[z]$,
     on which $M_2({\Bbb C}[z])$ acts by multiplication and
      $\partial_z$ acts by formal differentiation,  and
  $Y$ be the partially deformation-quantized space
   $Q_{\lambda}\mbox{\boldmath $\Omega$}_{{\Bbb A}^2/{\Bbb A}^1}$
   associated to the ring
   $S_{\lambda} :=
    {\Bbb C}\langle u,v,w \rangle/([v,w], [u,v], [u,w]-\lambda)$,
    where $\lambda\in{\Bbb C}$.
 Note\footnote{Also, we take the convention that
                $\partial_z\cdot m$ means the product in
                 ${\Bbb C}\langle M_2({\Bbb C}[z]), \partial_z\rangle$
                 and
                $\partial_z m$ means entry-wise formal differentiation
                 of $m$, for $m\in M_2({\Bbb C}[z])$.}
  that the action of $\partial_z$ on $N$
  induces an action of $\partial_z$ on $M_2({\Bbb C}[z])$
  by the entry-wise formal differentiation and
 the ${\Bbb A}^2/{\Bbb A}^1$ corresponds to
  ${\Bbb C}[v]\hookrightarrow{\Bbb C}[v,w]$.
 Consider the following special class of morphisms:
  $$
   \hspace{10em}
   \begin{array}{ccc}
    X & \xymatrix{\ar[rrr]^-{\varphi_{(A,B)}} &&&}          & Y \\[.6ex]
    R & \xymatrix{&&& \ar[lll]_-{\varphi_{(A,B)}^{\sharp}}}
      & S_{\lambda}\\[.6ex]
    \lambda\partial_z + A
      & \xymatrix{&&& \ar @{|->}[lll]} & u \\[.6ex]
    B & \xymatrix{&&& \ar @{|->}[lll]} & v \\[.6ex]
    z & \xymatrix{&&& \ar @{|->}[lll]} & w \\[.6ex]
   \end{array}\,,\hspace{2em}
   \mbox{$A,\, B\;\in\; M_2({\Bbb C}[z])\,$,}
  $$
  subject to
  $[\lambda\partial_z+A, B]\,=\,0\,$.
  (The other two constraints,
   $[B,z]\,=\,0\,$ and $\,[\lambda\partial_z+A, z]-\lambda\,=\,0\,$,
  are automatic.)
 Let
  $$
   A\; =\;
    \left[\begin{array}{cc} a_1 & a_2 \\ a_3 & a_4 \end{array}\right]
   \hspace{2em}\mbox{and}\hspace{2em}
   B\; =\;
    \left[\begin{array}{cc} b_1 & b_2 \\ b_3 & b_4 \end{array}\right]\,,
  $$
  where $a_i$, $b_j\in {\Bbb C}[z]$ and assume that $\lambda\ne 0$.
 Then, the associated system
   $\lambda\partial_zB+[A,B]=0$
   of homogeneous linear ordinary differential equations on $B$
   has a solution
  if and only if $A$ satisfies
   $$
    (a_1-a_4)^2 + 4 a_2a_3\;=\;0\,.
   $$
  Under this condition on $A$, the system has four fundamental solutions:
   $$
    \begin{array}{lcl}
    B_1  & =
     & \left[ \begin{array}{llll}
             1+\lambda^{-2}a_2a_3z^2
              &&& \lambda^{-1}a_2z
                  -\frac{1}{2}\lambda^{-2}(a_1-a_4)a_2 z^2
                  \hspace{2.3ex}  \\[.6ex]
             -\lambda^{-1}a_3z -\frac{1}{2}\lambda^{-2}(a_1-a_4)a_3 z^2
              &&& -\lambda^{-2}a_2a_3 z^2
            \end{array}
       \right]\,, \\[4ex]
    B_2  & =
     & \left[ \begin{array}{lll}
             \lambda^{-1}a_3z - \frac{1}{2}\lambda^{-2}(a_1-a_4)a_3 z^2
              && \hspace{1.6em}
                 1 - \lambda^{-1}(a_1-a_4)z -\lambda^{-2}a_2a_3 z^2\\[.6ex]
             - \lambda^{-2}a_3^2 z^2
               && \hspace{1.6em}
                  -\lambda^{-1}a_3z +\frac{1}{2}\lambda^{-2}(a_1-a_4)a_3z^2
            \end{array}
       \right]\,, \\[4ex]
    B_3  & =
     & \left[ \begin{array}{lll}
             -\lambda^{-1}a_2z -\frac{1}{2}\lambda^{-2}(a_1-a_4)a_2 z^2
              && \hspace{1.3ex}
                 - \lambda^{-2}a_2^2 z^2 \\[.6ex]
             1 + \lambda^{-1}(a_1-a_4)z - \lambda^{-2}a_2a_3 z^2
              && \hspace{1.3ex}
                 \lambda^{-1}a_2z + \frac{1}{2}\lambda^{-2}(a_1-a_4)a_2z^2
                 \hspace{1em}
            \end{array}
       \right]\,, \\[4ex]
    B_4  & =
     & \left[ \begin{array}{llll}
             -\lambda^{-2}a_2a_3 z^2
              &&& \hspace{1.3ex}
                  - \lambda^{-1}a_2z
                  + \frac{1}{2}\lambda^{-2}(a_1-a_4)a_2 z^2
                  \hspace{.9ex} \\[.6ex]
             \lambda^{-1}a_3z + \frac{1}{2}\lambda^{-2}(a_1-a_4)a_3 z^2
              &&& \hspace{1.3ex}
                  1 + \lambda^{-2}a_2a_3 z^2
            \end{array}
       \right]\,.
    \end{array}
   $$
  Denote this solution space by ${\Bbb C}^4_A$ with coordinates
   $(\hat{b}_1,\,\hat{b}_2,\,\hat{b}_3,\,\hat{b}_4)$
   and the correspondence
   $$
    (\hat{b}_1,\,\hat{b}_2,\,\hat{b}_3,\,\hat{b}_4)\;\;\;
     \longleftrightarrow\;\;\;
     \hat{b}_1 B_1\,+\,\hat{b}_2B_2\,+\,\hat{b}_3B_3\,+\,\hat{b}_4B_4\;
     =:\; B_{(\hat{b}_1,\hat{b}_2,\hat{b}_3,\hat{b}_4)}\,.
   $$
  Then,
   \begin{itemize}
    \item[$\cdot$] {\it
     the degree-$0$ term
     $B_{(0)}$ of
     $B=B_{(\hat{b}_1,\hat{b}_2,\hat{b}_3,\hat{b}_4)}$ (in $z$-powers)
     is given by
     $\left[\begin{array}{ll}
             \hat{b}_1 & \hat{b}_2 \\[.6ex] \hat{b}_3 & \hat{b}_4
            \end{array} \right]\,$,}

    \item[$\cdot$] {\it
     the characteristic polynomial of $B$ is identical to
     that of $B_{(0)}$.}
   \end{itemize}
 It follows that the image $\Image\varphi_{(A,B)}$ of $\varphi_{(A,B)}$
  is a (complex-)codimension-$1$ sub-quantum scheme in $Y$
  whose associated ideal in $S_{\lambda}$ contains the ideal
   $$
    \left( v^2-\trace B_{(0)}\,v + \determinant B_{(0)} \right)\,.
   $$

 Let $\mu_-$ and $\mu_+$ be the eigen-values of $B_{(0)}$.

 \bigskip

 \noindent
 {\it Case $(a):\,$ $\nu_-\ne\nu_+$.}\hspace{1ex}
  In this case,  the above ideal $((v-\nu_-)(v-\nu_+))$
   coincides with $\Ker\varphi_{(A,B)}^{\sharp}$ and, hence,
   describes precisely $\Image\varphi_{(A,B)}\subset Y$.
  Since $\varphi_{(A,B)}^{\sharp}(v)=B$,
   let $N_-:= \Ker(B-\nu_-)\subset N$.
  This is a rank-$1$ ${\Bbb C}[z]$-submodule of
    ${\Bbb C}[z]\oplus {\Bbb C}[z]$
    that is invariant also under
     $\varphi_{(A,B)}^{\sharp}(S_{\lambda})$.
  This gives $N_-$ a $S_{\lambda}/(v-\nu_-)$-module structure
   that has rank-$1$ as ${\Bbb C}[w]$-module.
  Similarly,
   $N_+:= \Ker(B-\nu_+)\subset N$ is invariant under
     $\varphi_{(A,B)}^{\sharp}(S_{\lambda})$
    and has a $\varphi_{(A,B)}^{\sharp}$-induced
     $S_{\lambda}/(v-\nu_+)$-module structure that is of rank-$1$
      as ${\Bbb C}[w]$-module.
  Let
   $$
    \begin{array}{crl}
     Z & :=
       & \Image\varphi_{(A,B)}\;\;
         =\;\; \Space (S_{\lambda}/((v-\nu_-)(v-\nu_+)))\\[.6ex]
       & =
       & \Space (S_{\lambda}/(v-\nu_-))
            \cup \Space (S_{\lambda}/(v-\nu_+))\;\;
         =:\;\; Z_-\cup Z_+
    \end{array}
   $$
   be the two connected components of the quantum subscheme
   $\Image\varphi_{(A,B)}\subset Y$  and
  denote the ${\cal O}_{Z_-}$-modules associated to $N_-$ and $N_+$
   by $(_{S_{\lambda}}N_-)^{\sim}$ and $(_{S_{\lambda}}N_+)^{\sim}$
   respectively.
  Then
   $$
    \varphi_{(A,B),\ast}{\cal E}\;
     =\; (_{S_{\lambda}}N_-)^{\sim}\oplus (_{S_{\lambda}}N_+)^{\sim}
    \hspace{1em}\mbox{with
     $(_{S_{\lambda}}N_-)^{\sim}$ supported on $Z_-$ and
     $(_{S_{\lambda}}N_+)^{\sim}$ on $Z_+$}\,.
   $$

 \bigskip

 \noindent
 {\it Case $(b):\,$ $\nu_-=\nu_+=\nu$.}\hspace{1ex}
 In this case, $\Ker\varphi_{(A,B)}^{\sharp}$ can be either
  $(v-\nu)$ or $((v-\nu)^2)$ and both situations happen.
 \begin{itemize}
  \item[$\cdot$]
   When $\Ker\varphi_{(A,B)}^{\sharp} = (v-\nu)$,
   $N={\Bbb C}[z]\oplus{\Bbb C}[z]$
    has a $\varphi_{(A,B)}^{\sharp}$-induced
    $S_{\lambda}/(v-\nu)$-module structure  and
   $\varphi_{(A,B),\ast}{\cal E}$ has support
    $\Image\varphi_{(A,B)}=\Space(S_{\lambda}/(v-\nu))\subset Y$.

  \item[$\cdot$]
   When $\Ker\varphi_{(A,B)}^{\sharp} = ((v-\nu)^2)$,
   $N={\Bbb C}[z]\oplus{\Bbb C}[z]$
    has a $\varphi_{(A,B)}^{\sharp}$-induced
    $S_{\lambda}/((v-\nu)^2)$-module structure  and
   $\varphi_{(A,B),\ast}{\cal E}$ has support
    $Z:=\Image\varphi_{(A,B)}=\Space(S_{\lambda}/((v-\nu)^2))\subset Y$.
   It contains an ${\cal O}_Z$-submodule
    $(_{S_{\lambda}}N_0)^{\sim}$, associated to
    $N_0:=\Ker(v-\nu)\subset N$, that is supported on
    $Z_0:=\Space(S_{\lambda}/(v-\nu))\subset Z$.
   In other words, in the current situation, $\varphi_{(A,B),\ast}{\cal E}$
    not only is of rank-$2$ as a ${\Bbb C}[w]$-module
    but also has a built-in $\varphi_{(A,B)}$-induced filtration
    $(_{S_{\lambda}}N_0)^{\sim}\subset \varphi_{(A,B),\ast}{\cal E}$.
 \end{itemize}

 \bigskip

 \noindent
 Thus, by varying $(A,B)$ in the solution space of
   $\lambda\partial_zB+[A,B]=0$
  so that the eigen-values of $B_{(0)}$ change from being distinct
   to being identical and vice versa,
 one realizes the Higgsing and un-Higgsing phenomena of D-branes
   in superstring theory for the current situation
  as deformations of morphisms from Azumaya quantum schemes to
   the open-string quantum target-space $Y$:

 
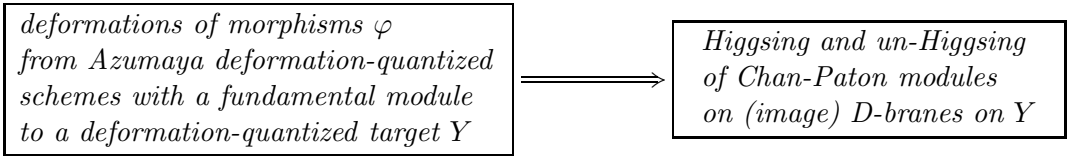
\bigskip
 \bigskip

 \hspace{1ex}
 \xymatrix{
  \framebox[17.6em][c]{\parbox{16.6em}{\it
   deformations of morphisms $\varphi$\\
   from Azumaya deformation-quantized\\ schemes
   with a fundamental module\\
   to a deformation-quantized target $Y$}}
   \ar @2{->}[rr]
   && \framebox[13.6em][c]{\parbox{11.6em}{\it
       Higgsing and un-Higgsing\\
       of Chan-Paton modules\\
       on (image) D-branes on $Y$}}
 } 

 \bigskip
 \bigskip

 \noindent
 Cf.~[L-L-S-Y: {\sc Figure~2-1-1}] for a similar phenomenon.

 This concludes the example.
}\end{example}


\bigskip

\subsection{Deformation quantizations of spectral covers
            in a cotangent bundle.}
We employ the notions from the previous subsection
 to discuss the notion of ``quantum spectral covers"\footnote{The
                         current subsection is written with
                          the particular works [D-H-S-V] and [D-H-S] of
                          Dijkgraaf, Hollands, Su{\l}kowski, and Vafa
                          in mind.
                         We thank Cumrun Vafa for the illuminations
                          of [D-H-S-V].
                         These works involve several mathematical themes.
                         Here we focus on a particular one:
                          the notion of {\it quantum spectral curves
                           from the viewpoint of D-branes}.
                         For that reason, it is not very appropriate to
                          attach a sub-title like
                          {\it Dijkgraaf-Holland-Su{\l}kowski-Vafa vs.\
                               Polchinski-Grothendieck}
                          to this subsection
                          though this is indeed what this subsection
                           is meant to be for the relevant part of
                           [D-H-S-V] and [D-H-S].
                         Readers are referred ibidem and references therein
                          for related stringy contents/pictures.}
 from the viewpoint of Azumaya geometry and
 the Polchinski-Grothendieck Ansatz.
A special case of this gives the notion of deformation quantizations
 of spectral curves in Hitchin's integrable systems.

\bigskip

\begin{flushleft}
{\bf A {\boldmath $1$}-parameter family of deformation quantizations
     of the cotangent bundle {\boldmath $\Omega$}$_W$.}
\end{flushleft}
Let
 $W$ be a smooth variety of dimension $n$ over ${\Bbb C}$,
 $\Omega_W$ be its sheaf of K\"{a}hler differentials,  and
 {\boldmath $\Omega$}$_W$
   be the total space
    $\boldSpec(\Sym^{\bullet\,}\Omega_W^{\vee})
                            = \boldSpec(\Sym^{\bullet\,}\Theta_W)$
    of $\Omega_W$.
One may construct
 a {\it $1$-parameter family of deformation quantizations}
 of {\boldmath $\Omega$}$_W$ as follows.

Let $p\in W$ be a geometric point on $W$.
Then there exists a Zariski open neighborhood $U$ of $p$ in $W$
 such that
  $\Omega_U$ is a free ${\cal O}_U$-module and that
  there admits an \'{e}tale morphism
  $\pi:U\rightarrow {\Bbb A}^n=\Spec({\Bbb C}[w_1,\,\cdots\,,\,w_n])$.
Denote the lifting of $w_i$ and $\partial_{w_i}$ on ${\Bbb A}^n$
 to $U$ under $\pi$ also by $w_i$ and $\partial_{w_i}$ respectively,
 for $i=1,\,\ldots\,,\, n$.
Then as both $U$ and $\pi$ are smooth,
 ${\cal D}_U$ as an ${\cal O}_W(U)$-algebra given by
 ${\Bbb C}\langle{\cal O}_W(U),\,
                  \partial_{w_1},\,\cdots\,,\,\partial_{w_n}\rangle$.
This is abstractly the algebra
 ${\Bbb C}\langle{\cal O}_W(U),\, p_1,\,\cdots\,,\,p_n\,\rangle/I$
 with $I$ the two-sided ideal
  $([p_i,p_j],\,[p_i,w_j]-\delta_{ij}\,:\, 1\le i, j\le n)$.
Here, we think of $p_i$ as a local section of the tangent sheaf
 ${\cal T}_W$ of $W$
 without a pre-assigned action on ${\cal O}_W$.
Note that,
 as ${\cal O}_W(U)$ is integral over an open subset of ${\Bbb A}^n$
     under $\pi$  and
    both $U$ and $\pi$ are smooth,
 the set of equations $[p_i,w_j]=\delta_{ij}$, $1\le i, j\, \le n$,
  determine the commutator $[p_i, f]\in {\cal O}_W(U)$,
  which is $\partial_{w_i}f$,
  for all $f\in {\cal O}_W(U)$ and $i=1,\,\ldots\,,\,n\,$.

\smallskip

\begin{notation}
{\bf [unital associative algebra generated by module].} {\rm
 (1)
  Let
  $S$ be a commutative ring and
  $R$ be a commutative $S$-algebra with a built-in $S\subset R$, and
  $M$ be a finitely generated $R$-module.
 Denote by $S\langle M\rangle$
  the unital associative $S$-algebra generated by elements
   of $M$ with the requirement that $S$ be in the center, and
  by $S\langle R,M\rangle$ be the unital associative $S$-algebra
   generated by $R\cup S\langle M\rangle$
   with a built-in $S$-algebra inclusions,
    $R\subset S\langle R,M\rangle$ and
    $S\langle M\rangle \subset S\langle R,M\rangle$.
 Note that $S$ is in the center of $S\langle R,M\rangle$
  while $R$ in general is not.

 (2)
 Let
  $Z$ be a scheme over a base ${\Bbb C}$-scheme $B$  and
  ${\cal F}$ be a coherent ${\cal O}_Z$-module.
 Denote by ${\cal O}_B\langle{\cal O}_Z,\,{\cal F}\rangle$
  the sheaf of unital associative ${\cal O}_Z$-algebras
  from the enlargement of ${\cal O}_Z$ by elements of ${\cal F}$
  with the requirement that the built-in
   ${\cal O}_B\subset {\cal O}_B\langle{\cal O}_Z,\,{\cal F}\rangle$
   be in the center.
 Over an affine open subset $U$ of $Z$
   that sits over an affine open subset $V$ of $B$,
  ${\cal O}_B\langle{\cal O}_Z,\,{\cal F}\rangle(U)$
  is the unital associative algebra
  ${\cal O}_B(V)\langle {\cal O}_Z(U), {\cal F}(U) \rangle$.
 Note that the image of the built-in inclusion
  ${\cal O}_Z\subset {\cal O}_B\langle {\cal O}_Z,\,{\cal F}\rangle$
  in general does not lie in the center.

 (3)
 Let
  $\Theta_{Z/B}$ be the sheaf of ${\cal O}_B$-derivations on $Z/B$  and
  ${\cal T}_{Z/B}$ be the relative tangent sheaf of $Z/B$.
 They are canonically isomorphic ${\cal O}_Z$-modules.
 However, for convenience, we take the convention that
  ${\cal O}_B\langle{\cal O}_Z,\,\Theta_{Z/B}\rangle$
   is the sheaf ${\cal D}_{Z/B}$ of algebras of
   differential operators on ${\cal O}_{Z/B}$
  (i.e.\  the $\Theta_{Z/B}$-action on ${\cal O}_Z$ via derivations
   is already included into its construction
   by setting $[p_i, f]=\partial_{w_i}f$)
  and that
   ${\cal O}_B\langle{\cal O}_Z,\,{\cal T}_{Z/B}\rangle$
   is constructed as Item (2) above,
   with ${\cal T}_{Z/B}$ treated only as an abstract coherent
   ${\cal O}_ Z$-module.
}\end{notation}

\smallskip

Let ${\cal T}_{({\Bbb A}^1\times W)/{\Bbb A}^1}$
 be the relative tangent sheaf of $({\Bbb A}^1\times W)/{\Bbb A}^1$.
Consider the ${\cal O}_{{\Bbb A}^1\times W}$-algebra
 ${\cal O}_{{\Bbb A}^1}
     \langle {\cal O}_{{\Bbb A}^1\times W},\,
             {\cal T}_{({\Bbb A}^1\times W)/{\Bbb A}^1}\rangle$.
Here, we take ${\Bbb A}^1$ as $\Spec({\Bbb C}[\lambda])$.
Let $I$ be the two-sided ideal sheaf of
 ${\cal O}_{{\Bbb A}^1}
     \langle {\cal O}_{{\Bbb A}^1\times W},\,
            {\cal T}_{({\Bbb A}^1\times W)/{\Bbb A}^1}\rangle$
 whose value over ${\Bbb A}^1\times U$,
  for $U$ being an affine open subset of $W$
   over which ${\cal T}_W$ is trivialized by
   $(p_1,\,\cdots\,,\,p_n)$ corresponding to
   $(\partial_{w_1},\,\cdots\,,\, \partial_{w_n})$,
 is given by
 $$
  ([p_i,p_j],\, [p_i,w_j]-\lambda\,:\, 1\le i,j\le n)\,.
 $$
Note that, by definition, $[\lambda,p_i]=0$, for $i=1,\,\ldots\,,\,n$,
  and that
 $[p_i,f] = \lambda\,\partial_{w_i}f$,
   where $f\in {\cal O}_{{\Bbb A}^1\times W}({\Bbb A}^1\times U)$,
  in the quotient ${\cal O}_{{\Bbb A}^1\times W}$-algebra
         ${\cal O}_{{\Bbb A}^1}
           \langle {\cal O}_{{\Bbb A}^1\times W},\,
                   {\cal T}_{({\Bbb A}^1\times W)/{\Bbb A}^1}\rangle
             /I$.
This gives a noncommutative space
 $$
  Q_{{\Bbb A}^1}\mbox{\boldmath $\Omega$}_W\;
   :=\; \Space(
         {\cal O}_{{\Bbb A}^1}
           \langle {\cal O}_{{\Bbb A}^1\times W},\,
                   {\cal T}_{({\Bbb A}^1\times W)/{\Bbb A}^1}\rangle
             /I )
 $$
 over ${\Bbb A}^1$.
It has the following properties:
 \begin{itemize}
  \item[(1)]
   The fiber $Q_{\lambda}\mbox{\boldmath $\Omega$}_W$
    over $\lambda\ne 0$ is isomorphic to the noncommutative space
    $\Space {\cal D}_W$.

 \item[(2)]
  The fiber over $\lambda=0$ is the commutative scheme
  {\boldmath $\Omega$}$_W$.

 \item[(3)]
  There is no local section of
   ${\cal O}_{{\Bbb A}^1}
     \langle {\cal O}_{{\Bbb A}^1\times W},\,
             {\cal T}_{({\Bbb A}^1\times W)/{\Bbb A}^1}\rangle$
   that is annihilated, either from the left or from the right,
    by a non-zero element of ${\Bbb C}[\lambda]$.
  Thus, we may think of
   $Q_{{\Bbb A}^1}\mbox{\boldmath $\Omega$}_W/{\Bbb A}^1$
   as a flat family of generically noncommutative spaces
    over ${\Bbb A}^1$, parameterized by $\lambda$.
\end{itemize}

\smallskip

\begin{definition}
{\bf [canonical family of deformation quantization].} {\rm
 We shall call the noncommutative space
  $Q_{{\Bbb A}^1}\mbox{\boldmath $\Omega$}_W$ over ${\Bbb A}^1$
  the {\it canonical family of deformation quantizations}
  of $\mbox{\boldmath $\Omega$}_W$.
}\end{definition}

\bigskip

\begin{flushleft}
{\bf Spectral covers via
     fibered morphisms from Azumaya schemes.}\footnote{Readers
                      are referred to [Hi1], [B-N-R], [Ox], [Ni],
                       and [Ma], [Don], [Do-M]
                       for the classical study of {\it Higgs pairs} and
                       their associated {\it spectral curves/covers}.
                      The current theme continues the discussion of
                       the theme
                       ``{\it Comparison with the spectral cover
                              construction and the Hitchin system}"
                       in [L-Y1: Sec.~4.1].
                      Here, we see one more example of the ubiquity
                       of Azumaya geometry in mathematics and
                       its recovering of D-brane
                       phenomena.
                      As illustrated in the precedent D(1) - D(4),
                      Azumaya geometry is a very fundamental
                       nature and geometry
                       a D-brane world-volume carries.
                      It gives the common origin of
                       many of the D-brany phenomena.
                      Furthermore, like what happens here,
                      such a structure is actually hidden
                       in many mathematical problems as well.
                      Despite the introduction of the notion of
                       ``maximally central algebra"
                        by Prof.\ Goro Azumaya in [Az] in year 1951,
                        which later came to be called ``Azumaya algebra",
                         and the study of it from the viewpoint
                         of algebras and representation theory,
                       the investigation of it as a geometric object
                        started only much later,
                        cf.\ related reference in
                        [L-Y1], [L-Y2], and [L-Y3].
                      The full richness of Azumaya geometry remains
                       to be explored.}$^{,}\,$\footnote{From
                      C.-H.L.$\;$:
                      The setting here rewrites and generalizes
                       some discussions
                       with {\it Mihnea Popa}
                       on the connection between
                        the D-brane Higgsing/un-Higgsing phenomenon
                        and spectral covers in spring 2002
                        ([Liu] and [Popa]).
                      During years 2001 - 2005,
                       Mihnea was giving lectures on a wide span of topics
                        in algebraic geometry:
                        from Grothendieck's foundation to
                        geometric topics at the frontier,
                       while {\it Shiraz Minwalla}
                        was giving lectures on an equally wide span
                        of topics in theoretical high energy physics:
                         from quantum field theory and
                          supersymmetry foundation
                         to stringy topics at the frontier.
                      Their systematic lectures in these four years
                        were filled with insight and enthusiasm  and
                        played a definite role
                       for the revival of the project in early 2007.}
\end{flushleft}
Let
 $W$ be as above and
 ${\cal E}$ and ${\cal N}$ be coherent locally free ${\cal O}_W$-modules.
Then one has the canonical isomorphisms
 $$
  \Hom_{{\cal O}_W}({\cal E},{\cal E}\otimes{\cal N})\;
   \simeq\; \Gamma({\cal E}^{\vee}\otimes{\cal E}\otimes{\cal N})\;
   \simeq\; \Hom_{{\cal O}_W}
             ({\cal N}^{\vee},\,\Endsheaf_{{\cal O}_W}({\cal E}))\,.
 $$
Here, $(\,\bullet\,)^{\vee}$ denotes the dual ${\cal O}_W$-module of
 $(\,\bullet\,)$.

\smallskip

\begin{definition}
{\bf [commutativity-admissible
      $\phi:{\cal E}\rightarrow {\cal E}\otimes{\cal N}$].}
{\rm
 An ${\cal O}_W$-module homomorphism
  $\phi:{\cal E}\rightarrow {\cal E}\otimes{\cal N}$
  is said to be {\it commutativity-admissible}
 if its corresponding
  $\widehat{\phi}:
    {\cal N}^{\vee}\rightarrow \Endsheaf_{{\cal O}_W}({\cal E})$
   has image contained in a commutative ${\cal O}_W$-subalgebra
   of $\Endsheaf_{{\cal O}_W}({\cal E})$.
}\end{definition}

\smallskip

A commutativity-admissible
 $\phi:{\cal E}\rightarrow {\cal E}\otimes{\cal N}$
 induces an ${\cal O}_W$-algebra homomorphism
 $$
  \varphi^{\sharp}\;:\;
   \Sym^{\bullet\,}({\cal N}^{\vee})\;\longrightarrow\;
    \Endsheaf_{{\cal O}_W}({\cal E})\,,
 $$
 which defines a morphism $\varphi$
  from the Azumaya scheme with a fundamental module
  $( W^{A\!z},
     {\cal E}):= (W, {\cal O}_W^{A\!z}:=\Endsheaf_{{\cal O}_W}({\cal E}),
     {\cal E})$
  to the total space
   {\boldmath ${\cal N}$}$=\boldSpec(\Sym^{\bullet\,}({\cal N}^{\vee}))$
   of the ${\cal O}_W$-module ${\cal N}$.
Note that, in this case,
 both $W^{A\!z}$ and {\boldmath ${\cal N}$} are spaces over $W$  and
 $\varphi:(W^{A\!z},{\cal E})\rightarrow$ {\boldmath ${\cal N}$}
  is a morphism between spaces over $W$.
Let $\pi:\,\mbox{\boldmath ${\cal N}$}\,\rightarrow W$ be the built-in
 morphism.
Then, by construction, $\pi_{\ast}\varphi_{\ast}{\cal E}\simeq{\cal E}$
 canonically.

\smallskip

\begin{lemma}
{\bf [$\phi$ vs.\ $\varphi$].} {\it
 Given coherent locally free ${\cal O}_W$-modules ${\cal E}$ and ${\cal N}$,
 there is a canonical one-to-one correspondence
  $$
   \left\{
    \begin{array}{l}
     \mbox{commutativity-admissible}\\
     \mbox{${\cal O}_W$-module homomorphisms}\\
     \phi:{\cal E}\rightarrow {\cal E}\otimes{\cal N}
    \end{array}
   \right\}\;
    \longleftrightarrow\;
   \left\{
    \begin{array}{l}
      \mbox{morphisms
       $\varphi:(W^{A\!z},{\cal E})\rightarrow$ {\boldmath ${\cal N}$}}\\
      \mbox{as spaces over $W$}
    \end{array}
   \right\}\,.
  $$
}\end{lemma}

\smallskip

\begin{lemma}
{\bf [generalized spectral cover from morphisms from Azumaya scheme].}
{\it
 When ${\cal N}$ is a line bundle ${\cal L}$ on $W$,
 any $\phi:{\cal E}\rightarrow {\cal E}\otimes{\cal L}$
  is commutativity-admissible  and
 the image scheme $\Image\varphi\subset$ {\boldmath ${\cal N}$}
  of the corresponding $\varphi$ lies in the spectral cover
  $\Sigma_{({\cal E},\phi)}$ in {\boldmath ${\cal N}$}
  associated to the pair $({\cal E},\phi)$.
 If furthermore $\Sigma_{({\cal E},\phi)}$ is smooth,
  then  $\Image\varphi=\Sigma_{({\cal E},\phi)}$.
}\end{lemma}

\smallskip

\noindent
This is consistent with the fact that, in the case of the Lemma,
 the morphism $\varphi$ carries the full information of
  the pair $({\cal E},\phi)$
 while the spectral cover $\Sigma_{({\cal E},\phi)}$ may not.
When ${\cal N}$ is a general coherent locally-free ${\cal O}_W$-module,
 then $\varphi$ can be thought of interchangeably
 as a {\it generalized Higgs pair} $({\cal E},\phi)$ for ${\cal N}$
  with $\phi$ commutativity-admissible.

\smallskip

\begin{remark}
{\it $[$general $\phi:{\cal E}\rightarrow {\cal E}\otimes{\cal N}$$]$.}
{\rm
 Let ${\cal E}$ and ${\cal N}$ be as above.
 Let ${\cal O}_W\langle{\cal N}^{\vee}\rangle$
  be the unital associative ${\cal O}_W$-algebra generated
  by ${\cal N}^{\vee}$
  with the requirement that ${\cal O}_W$ be in the center.
 Its associative ``space"
  $\Space({\cal O}_W\langle{\cal N}^{\vee}\rangle)$
  can be thought of as a noncommutative-affine-space bundle
  over $W$.
 A general ${\cal O}_W$-module homomorphism
  $\phi:{\cal E}\rightarrow {\cal E}\otimes{\cal N}$
  corresponds then to an ${\cal O}_W$-algebra homomorphism
  $\varphi^{\sharp}:{\cal O}_W\langle{\cal N}^{\vee}\rangle
       \rightarrow \Endsheaf_{{\cal O}_W}({\cal E})$.
 This can be thought of as defining a morphism
  $\varphi:(W^{A\,z},{\cal E}) \rightarrow
            \Space({\cal O}_W\langle{\cal N}^{\vee}\rangle)$
  over $W$.
}\end{remark}

\bigskip

\begin{flushleft}
{\bf Deformation quantization of spectral covers via morphisms
     from Azumaya schemes
     with a fundamental module with a flat connection.}
\end{flushleft}
Let
 $W$ be a smooth variety over ${\Bbb C}$,
 ${\cal E}$ and ${\cal N}$ be coherent locally free ${\cal O}_W$-modules,
 $\phi:{\cal E}\rightarrow {\cal E}\otimes\Omega_W$
  be a commutativity-admissible ${\cal O}_W$-module homomorphism, and
 $\varphi:(W^{A\!z},{\cal E})\rightarrow$ {\boldmath $\Omega$}$_W$
  be the associated morphism between spaces over $W$.
Let ${\cal E}_{{\Bbb A}^1}$
 be the locally-free ${\cal O}_{{\Bbb A}^1\times W}$-module
 on ${\Bbb A}^1\times W$ from the pull-back of ${\cal E}$
 under the projection map ${\Bbb A}^1\times W\rightarrow W$.
Here we take
 ${\Bbb A}^1$ as $\Spec({\Bbb C}[\lambda])$  and
 ${\cal E}_{{\Bbb A}^1}$ as a constant family of ${\cal O}_W$-modules
  over ${\Bbb A}^1$.
Denote $\Spec({\Bbb C}[\lambda,\lambda^{-1}])$
 by ${\Bbb A}^1-\{\mathbf 0\}$
 and the restriction of ${\cal E}_{{\Bbb A}^1}$
  to over ${\Bbb A}^1-\{\mathbf 0\}$ by
  ${\cal E}_{{\Bbb A}^1-\{\mathbf 0\}}$.
For convenience,
 we will denote a number in ${\Bbb C}$ also by $\lambda$.

\smallskip

\begin{definition}
{\bf [{\boldmath $\lambda$}-connection].}
{\rm ([Ari: Definition~2.1 and Example~2.2].)}
{\rm
 For $\lambda\in{\Bbb C}$,
 a {\it $\lambda$-connection} on ${\cal E}$ is a ${\Bbb C}$-linear map
  $\nabla:{\cal E}\rightarrow{\cal E}\otimes\Omega_W$ which satisfies
  the {\it $\lambda$-Leibniz rule}:
  $$
   \nabla(fs)\;=\; \lambda\cdot s\otimes df\, +\, f\nabla s
  $$
  for any $f\in {\cal O}_W$, $\, s\in{\cal E}$.
 Note that a $0$-connection on ${\cal E}$
  is an ${\cal O}_W$-module homomorphism  and,
 for $\lambda\ne 0$, $\nabla$ is a $\lambda$-connection
  if and only if $\lambda^{-1}\nabla$ is an (ordinary) connection.
 A $\lambda$-connection $\nabla$ on ${\cal E}$ is said to be {\it flat}
   if the connection $\lambda^{-1}\nabla$ is flat.
}\end{definition}

\smallskip

\noindent
The notion of $\lambda$-connection was introduced by Deligne;
it gives an interpolation between a Higgs field and a connection
 on ${\cal E}$,
cf.~[Sim].

\smallskip

\begin{definition}
{\bf [$\lambda$-connection deformation of $\phi$].}
{\rm
 An ${\cal O}_{{\Bbb A}^1-\{\mathbf 0\}}$-module homomorphism
  $$
  \nabla\;:\; {\cal E}_{{\Bbb A}^1-\{\mathbf 0\}}\; \longrightarrow\;
    {\cal E}_{{\Bbb A}^1-\{\mathbf 0\}}
      \otimes \Omega_{
       (({\Bbb A}^1-\{\mathbf 0\})\times W)/({\Bbb A}^1-\{\mathbf 0\})}
  $$
  that satisfies$\,$:
  \begin{itemize}
   \item[(1)]
    on each ${\cal E}_{\lambda}:= {\cal E}_{{\Bbb A}^1}|_{\lambda}$
     over a closed point, parameterized by $\lambda$,
     of ${\Bbb A}^1-\{\mathbf 0\}$,\\
    $\nabla^{\lambda}:=\nabla|_{\lambda}:
      {\cal E}_{\lambda}\rightarrow {\cal E}_{\lambda}\otimes \Omega_W$
    is a $\lambda$-connection on ${\cal E}_{\lambda}\,$,

   \item[(2)]
    $\nabla|_{\lambda=0}\;=\;\phi\,$,
  \end{itemize}
  is called a {\it $\lambda$-connection deformation} of $\phi$.
 If furthermore each $\nabla^{\lambda}$ is flat,
  then $\nabla$ is called a {\it flat $\lambda$-connection deformation}
  of $\phi$.
}\end{definition}

\smallskip

Let
 $X_{{\Bbb A}^1}
  = \Space( {\cal D}_{({\Bbb A}^1\times W)/{\Bbb A}^1} )$,
 $\, ({\cal E}_{{\Bbb A}^1},\,
     \nabla\mbox{ on ${\cal E}_{{\Bbb A}^1-\{\mathbf 0\}}$})$
  be as above,  and
 $Y_{{\Bbb A}^1}=Q_{{\Bbb A}^1}\mbox{\boldmath $\Omega$}_W$
  be the canonical family of deformation quantizations
  of {\boldmath $\Omega$}$_W$.
Let
 $$
  \begin{array}{cccl}
   X_{{\Bbb A}^1}
    & \xymatrix{ \ar[rr]^-{\varphi_{{\Bbb A}^1}} && }
    & Y_{{\Bbb A}^1} \\[1.2ex]
   {\cal O}_{{\Bbb A}^1}
    \langle {\cal O}_{{\Bbb A}^1\times W},\,
            \Theta_{({\Bbb A}^1\times W)/{\Bbb A}^1}   \rangle
    & \xymatrix{ && \ar[ll]_-{\varphi_{{\Bbb A}^1}^{\sharp}} }
    & {\cal R}       \\[.6ex]
   \lambda & \xymatrix{ && \ar @{|->} [ll] } & \lambda \\[.6ex]
   w_i     & \xymatrix{ && \ar @{|->} [ll] } & w_i     \\[.6ex]
   \lambda\partial_{w_i} & \xymatrix{ && \ar @{|->} [ll] } & p_i
    &.
  \end{array}
 $$
Here we adopt the notation in the construction of
 $Q_{{\Bbb A}^1}\mbox{\boldmath $\Omega$}_W$.
Then, $\varphi_{{\Bbb A}^1}$ is a morphism of spaces over ${\Bbb A}^1$
 with the following the properties:
 \begin{itemize}
  \item[$\cdot$]
   $\varphi_0:= \varphi_{{\Bbb A}^1}|_{\lambda=0}$
    is the composition
    $\,(W^{A\!z,{\cal D}},{\cal E}^{\nabla^{\prime}})
        \longrightarrow\, (W^{A\!z},{\cal E})\,
         \stackrel{\varphi}{\longrightarrow}\,
          \mbox{\boldmath $\Omega$}_W\,$;
    here,
     as $\lambda^{-1}\nabla$ does not extend over $\lambda=0$,
     we take $\nabla^{\prime}$ to be an arbitrary auxiliary
      flat connection on ${\cal E}$ to render ${\cal E}$
      a ${\cal D}_W$-module and
     let
      $(W^{A\!z,{\cal D}},{\cal E}^{\nabla^{\prime}})
                          \longrightarrow (W^{A\!z},{\cal E})$
      be the built-in dominant morphism;

  \item[$\cdot$]
   $\varphi_\lambda := \varphi_{{\Bbb A}^1}|_{\lambda}\,
     :\: (W^{A\!z,{\cal D}},{\cal E}^{\lambda^{-1}\nabla_{\lambda}})
         \longrightarrow Q_{\lambda}\mbox{\boldmath $\Omega$}_W$,
   for $\lambda\in {\Bbb A}^1-\{\mathbf 0\}\,$.
 \end{itemize}
Note that
 $\varphi_{0\,\ast}{\cal E}=\varphi_{\ast}{\cal E}$
  is a ${\cal O}_{\mbox{\scriptsize\boldmath $\Omega$}_W}$-module
  flat over $W$ with relative dimension $0$ and
   relative length $=\rank{\cal E}$.
 $\Supp({\cal E})=\Image\varphi$,
  which is identical to the spectral curves $\Sigma_{({\cal E},\phi)}$
  when $W$ is a smooth curve $C$ and $\Sigma_{({\cal E},\phi)}$
  is smooth.
On the other hand, for $\lambda\ne 0$,
 $\varphi_{\lambda\ast}{\cal E}$ has support
 $\Image\varphi_{\lambda}$,
 which is the whole $Q_{\lambda}\mbox{\boldmath $\Omega$}_W$.
The characteristic variety for
 the ${\cal O}_{Q_{\lambda}\mbox{\scriptsize\boldmath $\Omega$}_W}$-module
 $\varphi_{\lambda\ast}{\cal E}$
 is the zero-section of {\boldmath $\Omega$}$_W/W$.

\smallskip

\begin{remark}
{\it $[$existence/interpretation of quantum spectral covers$]$.}
{\rm
 Due to the fact that the Weyl algebras are simple algebras,
  the spectral curve $\Sigma_{({\cal E},\phi)}$
  in {\boldmath $\Omega$}$_C$ in general do not have a direct
  deformation quantization into $Q_{\lambda}\mbox{\boldmath $\Omega$}_C$
  by the ideal sheaf of $\Sigma_{({\cal E},\phi)}$
  in ${\cal O}_{\mbox{\scriptsize\boldmath $\Omega$}_C}$
 since this will only give
  ${\cal O}_{Q_{\lambda}\mbox{\scriptsize\boldmath $\Omega$}_C}$,
   which corresponds to the empty subspace of
   $Q_{\lambda}\mbox{\boldmath $\Omega$}_C$.
 The setting above replaces the notion of quantum spectral curves by
  quantum deformation $\varphi_{\lambda}$ of the morphism $\varphi$
  associated to the Higgs/spectral pair $({\cal E},\phi)$.
}\end{remark}

\newpage
\baselineskip 13pt

\end{document}